# RANK ONE CASE OF DWORK'S CONJECTURE

DAQING WAN

## 1. Introduction

In the higher rank paper [17], we reduced Dwork's conjecture from higher rank case over any smooth affine variety $X$ to the rank one case over the simplest affine space $\mathbf{A}^n$. In the present paper, we finish our proof by proving the rank one case of Dwork's conjecture over the affine space $\mathbf{A}^n$, which is called the key lemma in [17]. The key lemma had already been proved in [16] in the special case when the Frobenius lifting $\sigma$ is the simplest $q$-th power map $\sigma(x) = x^q$. Thus, the aim of the present paper is to treat the general Frobenius lifting case. Our method here is an improvement of the limiting method in [16]. It allows us to move one step further and obtain some explicit information about the zeros and poles of the unit root L-function. As in [16], to handle the rank one case, we are forced to work in the more difficult infinite rank setting, see section 2 for precise definitions of the various basic infinite rank notions. Let $\mathbf{F}_q$ denote the finite field of characteristic $p > 0$. Our main result of this paper is the following theorem.

**Theorem 1.1.** *Let $\phi$ be a nuclear overconvergent $\sigma$-module over the affine $n$-space $\mathbf{A}^n/\mathbf{F}_q$, ordinary at the slope zero side. Let $\phi_{\mathrm{unit}}$ be the unit root (slope zero) part of $\phi$. Assume that $\phi_{\mathrm{unit}}$ has rank one. Let $\psi$ be another nuclear overconvergent $\sigma$-module over $\mathbf{A}^n/\mathbf{F}_q$. Then for each integer $k$, the L-function $L(\psi \otimes \phi_{\mathrm{unit}}^{\otimes k}, T)$ is $p$-adic meromorphic. Furthermore, the family $L(\psi \otimes \phi_{\mathrm{unit}}^{\otimes k}, T)$ of L-functions parametrized by integers $k$ in each residue class modulo $(q-1)$ is a strong family of meromorphic functions with respect to the $p$-adic topology of $k$.*

A finite rank $\sigma$-module is automatically nuclear. Thus, Theorem 1.1 includes the key lemma of [17] over $\mathbf{A}^n$ as a special case. The basic ideas of the present paper are the same as the limiting approach in [16]. The details are, however, quite different. In the simplest $q$-th power Frobenius lifting case, one has the fundamental Dwork trace formula available, which is completely explicit for uniform estimates. This makes it easy to extend the Dwork trace formula to infinite rank setting. It also makes it possible to see the various analytic subtleties involved in a concrete case. As a result, we were able to prove analytically optimal results in [16]. For a general Frobenius lifting (even over the simplest affine $n$-space $\mathbf{A}^n$ as we shall work in this paper), one has to use the much more difficult Monsky trace formula which is a generalization of Dwork's trace formula. Thus, the first task of this paper is to extend the Monsky trace formula to infinite rank setting and to make it sufficiently

1991 *Mathematics Subject Classification.* Primary 11G40, 11S40; Secondary 11M41, 14G15.
*Key words and phrases.* L-functions, Fredholm determinants, $p$-adic meromorphic continuation, nuclear $\sigma$-modules and Banach modules.
This work was partially supported by NSF.

explicit so that our previous explicit limiting approach can be carried out. As indicated in [16], there are two versions of the limiting approach. The first version is to stay in finite rank setting but to take the limit of an essentially continuous (not continuous) family of functions. This is the version adapted in [16]. The second version, briefly outlined in [16] but without giving a detailed proof there, is to work directly with infinite rank setting and continuous family of functions. Here, we shall follow the second version of the limiting approach and work out all the details. This version gives more transparent and conceptual results, such as an explicit formula (see Theorem 1.9 below) for the unit root L-function $L(\psi \otimes \phi_{\mathrm{unit}}^{\otimes k}, T)$ in terms of nuclear overconvergent $\sigma$-modules. This explicit formula is important to gain further information about the distribution of the zeros and poles of the unit root L-function $L(\psi \otimes \phi_{\mathrm{unit}}^{\otimes k}, T)$. In this direction, we have the following result.

**Theorem 1.2.** *Let $\phi$ be an overconvergent $\sigma$-module of some finite rank $r$ over the affine $n$-space $\mathbf{A}^n/\mathbf{F}_q$, ordinary at the slope zero side. Let $\phi_{\mathrm{unit}}$ be the unit root (slope zero) part of $\phi$. Assume that $\phi_{\mathrm{unit}}$ has rank one. Let $\psi$ be another finite rank overconvergent $\sigma$-module over $\mathbf{A}^n/\mathbf{F}_q$. Then, for all integers $k$, we can write*

$$L(\psi \otimes \phi_{\mathrm{unit}}^{\otimes k}, T) = \frac{f_1(k, T)}{f_2(k, T)}, \ f_1(k, 0) = f_2(k, 0) = 1,$$

*where $f_1(k, T)$ (resp. $f_2(k, T)$) is a family of uniformly continuous and uniformly entire functions whose Newton polygons lies above the graph of the function*

$$Q(x) = c_1 x^{1 + \frac{1}{n+r-1}} - c_2 x$$

*for some positive constants $\{c_1, c_2\}$ independent of $k$.*

As an application of Theorem 1.2, we obtain explicit information about the degree of each pure slope piece of the L-function $L(\psi \otimes \phi_{\mathrm{unit}}^{\otimes k}, T)$. Recall that the degree (resp. the total degree) of a rational function means the degree of the numerator minus (resp. plus) the degree of the denominator. For a non-negative rational number $s$, let $d_s(k)$ (resp. $D_s(k)$) denote the degree (resp. the total degree) of the slope $s$ part of the meromorphic function $L(\psi \otimes \phi_{\mathrm{unit}}^{\otimes k}, T)$. Although $D_s(k)$ is a non-negative integer, the integer $d_s(k)$ can be negative. The $p$-adic Riemann hypothesis in this situation is then to understand the two degree functions $d_s(k)$ and $D_s(k)$ for each given integer $k$. Theorem 1.2 implies the following uniform upper bound for $d_s(k)$ and $D_s(k)$.

**Theorem 1.3.** *In the situation of Theorem 1.2, there is a positive constant $c$ such that for all rational numbers $s \geq 0$ and all integers $k$, the number of reciprocal zeros and reciprocal poles of $L(\psi \otimes \phi_{\mathrm{unit}}^{\otimes k}, T)$ with slopes at most $s$ is bounded by*

$$\sum_{t \in [0, s]} D_t(k) \leq c(s+1)^{n+r-1}.$$

*In particular, we have the uniform polynomial bound*

$$|d_s(k)| \leq D_s(k) \leq c(s+1)^{n+r-1}.$$

As another application of Theorem 1.2, we derive, with the help of a continuity result and the method of [15], the following information about the variation of the degree function $d_s(k)$ for fixed slope $s$ when $k$ varies.

**Theorem 1.4.** *In the situation of Theorem 1.2, there is a positive integer c such that whenever $k_1$ and $k_2$ are two integers satisfying the congruence*

$$k_1 \equiv k_2 \ (\mathrm{mod}(q-1)p^{c[s+1]^{n+r}}),$$

*we have the equality*

$$d_t(k_1) = d_t(k_2),$$

*for all $0 \le t \le s$.*

*Remark 1.5.* It is probably unreasonable to expect a similar result for the total degree function $D_s(k)$ due to possible cancellation of zeros and poles.

Theorem 1.4 generalizes a weak version of the Gouvêa-Mazur conjecture from the universal family of elliptic curves to an arbitrary family of algebraic varieties with one unit root fibre by fibre. An interesting higher dimensional geometric example is a good family of Calabi-Yau hypersurfaces where it is known that the middle dimensional relative crystalline cohomology (taken to be our $\phi$) has exactly one unit root fibre by fibre. In this case, it seems that the unit root L-function $L(\phi_{\mathrm{unit}}^{\otimes k}, T)$ can be globally constructed from the mirror map. It would be interesting to understand further connections between the two topics. Here, the crystal structure (the differential equation) of the relative crystalline cohomology would play an important role. Another interesting topic is to explore the possibility of using the special values $L(\phi_{\mathrm{unit}}^{\otimes k}, 1)$ parametrized by integer $k$ to construct the conjectural $p$-adic L-function attached to an algebraic variety defined over a number field. This relationship is not completely well understood even in the case of an elliptic curve, although some positive evidences are available. The bounds in Theorems 1.2-1.4 are probably asymptotically best possible in general. They can, however, be improved in some special cases, see the end of section 8 for a brief discussion of possible such improvements. In another direction, Coleman has shown that in the rank one case, the strong family assertion of $L(\psi \otimes \phi_{\mathrm{unit}}^k, T)$ in Theorem 1.1 can be significantly strengthened by extending the domain of the $p$-adic parameter $k$ and by using the theory of rigid analysis. This raises a host of new interesting questions, greatly generalizing the study of the eigencurve problem as initiated by Coleman-Mazur in [3] for the elliptic family.

If the unit root part $\phi_{\mathrm{unit}}$ is of rank greater than one, our result shows that the L-function $L(\psi \otimes \phi_{\mathrm{unit}}^{\otimes k}, T)$ is still meromorphic for each fixed integer $k$, but the corresponding statements in Theorems 1.2-4 are all false. In fact, in this higher rank situation, no non-trivial uniform bounds are possible due to the simple reason that the sequence $x^k$ as $k$ varies is unbounded for any fixed real number $x > 1$. The simplest counter-example is to take $\phi = \phi_{\mathrm{unit}}$ to be the direct sum of $r$ ($r > 1$) copies of the identity $\sigma$-module. If one replaces the $k$-th tensor power $\phi_{\mathrm{unit}}^{\otimes k}$ by the $k$-th iterate $\phi_{\mathrm{unit}}^k$, then there are a number of open questions in higher rank case, see section 8 in [17].

We now give a brief description of the content of the paper and the main steps in our proof. Section 2 contains a brief discussion of the foundational material and the infinite rank setup. It extends the easier finite rank setting of [17] to the more difficult infinite rank setting, improving the exposition in [16]. The concepts introduced are Banach modules over certain rings, topological basis, formal basis, orthonormal basis, $\sigma$-modules, nuclear $\sigma$-modules, overconvergent $\sigma$-modules, nuclear overconvergent $\sigma$-modules and the various tensor operations. The L-function is then defined for a nuclear $\sigma$-module. Section 3 discusses the uniform family version of overconvergent $\sigma$-modules and reviews several notions from [16] about a uniform family of entire (resp. meromorphic) functions.

Section 4 is devoted to Dwork operators and related explicit estimates. The main notion here is that of a contracting overconvergent Dwork operator. The main result of this section is the following theorem.

**Theorem 1.6.** *Let $\Theta$ be a contracting overconvergent Dwork operator. Then, its Fredholm determinant $det(I - T\Theta)$ is well defined and entire. More generally, for a uniform and continuous family $\Theta(k)$ of contracting overconvergent Dwork operators, the Fredholm determinant $det(I - T\Theta(k))$ is a uniform and continuous family of entire functions.*

In section 5, we extend the Monsky trace formula from finite rank case over $\mathbf{A}^n$ to infinite rank case over $\mathbf{A}^n$. The idea is to use the finite rank version as shown in [17] and to take a suitable limit. That is, we prove

**Theorem 1.7.** *Let $\phi$ be an overconvergent nuclear $\sigma$-module over $\mathbf{A}^n/\mathbf{F}_q$. Then, there are $n + 1$ contracting overconvergent Dwork operators $\Theta_i(\phi)$ such that*

$$(1.1) \qquad L(\phi, T) = \prod_{i=0}^{n} \det(I - T\Theta_i(\phi))^{(-1)^{i-1}}.$$

Putting Theorem 1.6 and Theorem 1.7 together, we deduce

**Theorem 1.8.** *Let $\phi$ be an overconvergent nuclear $\sigma$-module over $\mathbf{A}^n/\mathbf{F}_q$. Then, the L-function $L(\phi, T)$ is p-adic meromorphic. More generally, for a uniform and continuous family $\phi(k)$ of nuclear overconvergent $\sigma$-modules, the family of L-functions $L(\phi(k), T)$ parametrized by integers $k$ is a strong family of meromorphic functions with respect to the p-adic topology of $k$.*

In section 6, we review the definitions of the Newton polygon, basis polygon, ordinariness of nuclear $\sigma$-modules. In order to prove Theorem 1.1, it is sufficient to assume that $\phi$ is normalized in the sense that a suitable overconvergent matrix $B$ of $\phi$ is of the block form,

$$B = \begin{pmatrix} 1 + \pi b_{00} & \pi B_{01} \\ \pi B_{10} & \pi B_{11} \end{pmatrix},$$

where $\pi$ is the uniformizer and $1 + \pi b_{00}$ is an overconvergent $1 \times 1$ matrix.

In section 7, we begin the proof of Theorem 1.1. Without loss of generality, we assume that $\phi$ is normalized. Using the infinite rank version of the decomposition formula (1.5) in [17], we derive the following limiting formula:

$$(1.2) \qquad L(\psi \otimes \phi_{\mathrm{unit}}^k, T) = \lim_{m \to \infty} \prod_{i \geq 1} L(\psi \otimes \mathrm{Sym}^{k+p^m-i} \phi \otimes \wedge^i \phi, T)^{(-1)^{i-1} i}.$$

To prove that the limit exists and is meromorphic, the key is to construct a new overconvergent nuclear $\sigma$-module $(M_{\infty,k}, \phi_{\infty,k})$ called the limiting $\sigma$-module of the sequence $\mathrm{Sym}^{k+p^m} \phi$ $(m = 1, 2, \cdots)$ such that

$$(1.3) \qquad \lim_{m \to \infty} L(\psi \otimes \mathrm{Sym}^{k+p^m-i} \phi \otimes \wedge^i \phi, T) = L(\psi \otimes \phi_{\infty,k-i} \otimes \wedge^i \phi, T).$$

The tensor product of two nuclear overconvergent $\sigma$-modules is still nuclear and overconvergent. Thus, the right side of (1.3) is the L-function of a nuclear overconvergent $\sigma$-module and hence is meromorphic by Theorem 1.8. This implies that the unit root L-function $L(\psi \otimes \phi_{\text{unit}}^k, T)$ is also meromorphic. Furthermore, as a by-product, our proof gives the following stronger explicit formula.

**Theorem 1.9.** *Assume that we are in the situation of Theorem 1.1. Assume further that $\phi$ is normalized. Then, we have the following formula:*

$$(1.4) \qquad L(\psi \otimes \phi_{\text{unit}}^k, T) = \prod_{i \geq 1} L(\psi \otimes \phi_{\infty, k-i} \otimes \wedge^i \phi, T)^{(-1)^{i-1} i}.$$

Note that the product in (1.4) is essentially finite if we restrict to a finite disk. This is because $\phi$ is nuclear and thus $\wedge^i \phi$ is divisible by higher and higher power of the uniformizer $\pi$ as $i$ grows. Thus, there is a crucial shifting factor of $\pi$ in (1.4). The product in (1.4) is automatically finite if $\phi$ has finite rank since $\wedge^i \phi = 0$ for large $i$.

Using the explicit estimates in section 4, we derive in section 8, a uniform lower bound for the Newton polygon of the various Fredholm determinants arising from the family of L-functions $L(\psi \otimes \phi_{\text{unit}}^k, T)$ parametrized by $k$. This together with Theorem 1.9 then implies explicit information about the distribution of the zeros and poles of the L-function $L(\phi_{\text{unit}}^k, T)$ as stated in Theorems 1.2-4.

We would like to point out that the results of this paper do not include the results in [16] as a special case. In terms of the Frobenius lifting $\sigma$, this paper is more general and allows us to complete our proof of the rank case case of Dwork's conjecture over the affine $n$-space. However, from analytic point of view, the results of the present paper are not optimal and thus somewhat weaker than the results in [16]. The reason is that the ambient $\sigma$-module $\phi$ is only assumed to be $c$ log-convergent in [16] which is the weakest possible condition for the study of L-functions. In contrast, the present paper assumes that the ambient $\sigma$-module is overconvergent which is an analytically stronger condition, although it is enough for our immediate application to Dwork's conjecture. The overconvergent condition on $\phi$ simplifies some of the analytic arguments but our proof here is still quite analytic due to the nature of the problem.

**Acknowledgement**. The author would like to thank the referee for his careful reading and for sketching out some important foundational material on $p$-adic Banach modules which significantly clarified some of the basic concepts. Part of the material in this paper was worked out when the author was visiting the University of Rennes in the summer of 1998. He would like to thank the algebraic geometric group at Rennes for their warm hospitality.

## 2. Nuclear $\sigma$-modules and L-functions

Let $\mathbf{F}_q$ be the finite field of $q$ elements of characteristic $p$. Let $R$ be a complete discrete valuation ring with uniformizer $\pi$ and with residue field $\mathbf{F}_q$. Let $K$ be the quotient field of $R$. Thus, the ring $R$ consists of those elements $a \in K$ such that $\text{ord}_\pi a \geq 0$. We assume that $K$ has characteristic zero (the characteristic $p$ case is actually somewhat easier). One can then think of $K$ as a finite extension of the field $\mathbf{Q}_p$ of $p$-adic rational numbers with residue field $\mathbf{F}_q$. We normalize the $p$-adic ($\pi$-adic) absolute value on $K$ by defining $|\pi|_\pi = 1/p$. Let $n$ be a fixed positive integer. We shall consider various nuclear $\sigma$-modules over the affine $n$-space $\mathbf{A}^n$, which are certain Banach modules with a nuclear action of a semi-linear operator. Throughout the paper, our base space will be $\mathbf{A}^n$.

First, we define various coefficient rings. Let

$$A_0 = \{ \sum_{u \in \mathbf{Z}_{\geq 0}^n} a_u X^u | a_u \in R, \ \lim_{|u| \to \infty} |a_u|_\pi = 0 \}$$

be the ring of **convergent** power series over $R$, where for a lattice point $u = (u_1, \cdots, u_n) \in \mathbf{Z}_{\geq 0}^n$, we define $X^u = X_1^{u_1} \cdots X_n^{u_n}$ and $|u| = u_1 + \cdots + u_n$. The **overconvergent** subring $A$ of $A_0$ is defined by

$$A = \{ \sum_{u \in \mathbf{Z}_{\geq 0}^n} a_u X^u | a_u \in R, \ \lim_{|u| \to \infty} \inf \frac{\operatorname{ord}_\pi a_u}{|u|} > 0 \}.$$

Thus, the ring $A$ consists of those power series $f(X_1, \cdots, X_n)$ in $A_0$ which converge in the polydisk $|X_i|_\pi < r$ for some real number $r > 1$, where $r$ depends on $f$. Both $A_0$ and $A$ are Noetherian. Let $\sigma$ be a fixed $R$-algebra endomorphism of $A$ which lifts the $q$-th power Frobenius map of the polynomial ring $\mathbf{F}_q[X]$. Thus, there are unique elements $f_i \in A$ $(1 \leq i \leq n)$ such that

$$\sigma(X_i) = X_i^q + \pi f_i.$$

These relations also define an $R$-algebra endomorphism of $A_0$ by continuity. Although $\sigma$ acts on the bigger ring $A_0$, we remind the reader that our $\sigma$ is always originally defined over $A$. That is, the elements $\sigma(X_i)$ are always overconvergent.

The full ring $A_0$ is complete. This fact makes all algebraic definitions simpler. But it is too large in the sense that not much can be proved about the L-functions that we will define. The smaller overconvergent subring $A$ is important for us but it is not complete. This makes all algebraic definitions more complicated and one has to use analytic ideas. On the positive side, the ring $A$ is weakly complete in some sense. In nice situations, we can replace $A$ by a dense $R$-submodule which is not a ring any more but is complete. Such a complete $R$-submodule makes the analysis simpler. We will get to these submodules in later sections when dealing with finer properties. For now, we define our basic notions in terms of the big complete ring $A_0$. The ring $A_0$ serves as our ground ring or coefficient ring.

The ring $A_0$ is a complete normed $R$-algebra under the Gauss norm:

$$\| \sum_u a_u X^u \| = \max_u |a_u|_\pi.$$

Thus, the ring $A_0$ becomes a Banach algebra over $R$ and the monomials $\{X^u\}$ form an orthonormal basis of $A_0$ over $R$. We need to consider certain Banach modules over $A_0$ with certain nuclear $\sigma$-linear endomorphisms. A closely related theory of Banach modules over a ring (Banach algebra) with completely continuous maps is given in Coleman [2], generalizing the $p$-adic spectral theory of Serre [11] over the field $K$ and of Monsky [10] over the complete discrete valuation ring $R$. These theories, however, mostly restrict to Banach modules with an orthonormal basis and with linear maps, although some of the results may hold for more general Banach modules with certain semi-linear endomorphisms. The treatment below of the foundational material on Banach modules, influenced by the referee's comments, is an improvement of those given in [16]. Again, we only briefly touch on those properties which are relevant to our immediate applications. It may be interesting to give a full systematic development of the following foundational material which would be useful in some further study of L-functions.

**Definition 2.1.** An $A_0$-module $M$ is called a Banach module over $A_0$ if the following three properties are satisfied.

1. $M$ is a complete topological continuous $A_0$-module with a basis of open neighborhoods of the origin consisting of $A_0$-submodules.

2. $M$ has a norm $\|?\|$ satisfying the relations

$$\|m\| = 0 \text{ if and only if } m = 0, \ m \in M,$$

$$\|am\| = \|a\| \cdot \|m\|, \text{ for } a \in A_0, \ m \in M,$$

$$\|m_1 + m_2\| \leq \max(\|m_1\|, \|m_2\|), \text{ for } m_1, m_2 \in M.$$

3. If $\{m_i\}$ is a sequence in $M$ such that

$$\lim_i \|m_i\| = 0,$$

then the sequence $\{m_i\}$ converges to zero in $M$.

Note that in addition to the given topology on $M$, the norm on $M$ induces another topology on $M$, called the norm topology. The norm topology on $M$ may be different from the original topology on $M$. Condition (iii) in the above definition refers to a relationship between the two topologies on $M$. One can replace the ground ring $A_0$ by any Banach algebra and hence get a more general notion. Since, the ring $A_0$ is all we use in this paper, we will stick with this ground ring $A_0$. The simplest examples of Banach $A_0$-modules are $A_0$ itself and its quotient by a closed ideal. In particular, the ring $R$ is also a Banach $A_0$-module. We now examine two more examples.

**Example 2.2.** Let $I$ be a set. Let

$$b(I) = \{\{a_i\}_{i \in I} | a_i \in A_0\}$$

be the set of all sequences $\{a_i\}_{i \in I}$ of elements in $A_0$ indexed by $I$. If $I$ is the empty set, then $b(I)$ is defined to be the trivial zero module. Clearly, $b(I)$ is an $A_0$-module under componentwise addition and scalar multiplication. Define the norm on $b(I)$ by

$$\|a\| = \max_{i \in I} \|a_i\|.$$

Define the topology on $b(I)$ by requiring the submodules

$$M_S = \{\{a_i\}_{i \in I} : \|a_i\| < \epsilon \text{ for } i \in S\},$$

where $S$ is a finite subset of $I$ and $\epsilon$ is a positive real number, to be a basis of open neighborhoods of the origin. One checks that $b(I)$ is a Banach module over $A_0$. Furthermore, if $I$ is infinite, the topology on $b(I)$ is different from the norm topology on $b(I)$.

**Example 2.3.** Let $b(I)$ be the Banach $A_0$-module in the previous example. Let $c(I)$ be the $A_0$-submodule

$$c(I) = \{\{a_i\}_{i \in I} | a_i \in A_0, \lim_i \|a_i\| = 0\},$$

which is the set of all convergent sequences of elements in $A_0$ indexed by $I$. The norm $\|?\|$ on $c(I)$ is the same as the one on $b(I)$. But we define the topology on $c(I)$ by requiring the submodules

$$M_\epsilon = \{\{a_i\}_{i \in I} \in c(I) : \|a_i\| < \epsilon \text{ for } i \in I\},$$

where $\epsilon$ is a positive real number, to be a basis of open neighborhoods of the origin. In this way, the space $c(I)$ becomes a Banach $A_0$-module. One checks that the topology on $c(I)$ agrees with the norm topology on $c(I)$.

Let $M$ be a Banach module over $A_0$. A **topological basis** for $M$ is a subset $\vec{e} = \{\{e_i\}_{i \in I} | e_i \in M\}$ indexed by some set $I$ such that for every element $m \in M$, there exists a unique collection $\{a_i \in A_0\}_{i \in I}$ such that the net

$$\sum_{i \in S} a_i e_i$$

where $S$ ranges over finite subsets of $I$, converges to $m$ in $M$ (in term of the topology of $M$) and moreover,

$$\|m\| = \max_{i \in I} \|a_i\|.$$

In this case, we call the series

$$(2.1) \qquad\qquad m = \sum_{i \in I} a_i e_i$$

the expansion of $m$ with respect to $\vec{e}$. The cardinality of $I$ is called the rank of $M$. The rank is clearly independent of the choice of basis if it is finite. Otherwise, it is always infinite. In this paper, we always assume that the index set $I$ is at most countable. Thus, the rank is always independent of the choice of basis. The basis $\vec{e}$ is called a **formal basis** for $M$ if for every collection $\{a_i \in A_0\}_{i \in I}$, there exists a unique element $m \in M$ with the expansion (2.1). A typical such example is the module $b(I)$ given in Example 2.2. The basis $\vec{e}$ is called an **orthonormal basis** for $M$ if for every element $m \in M$, the expansion in (2.1) satisfies $\lim_i \|a_i\| = 0$. A typical such example is the module $c(I)$ given in Example 2.3.

For a finite rank free $A_0$-module $M$, all the above three notions on basis are the same concept. Any free $A_0$-module $M$ of finite rank can be viewed as a Banach module over $A_0$ of finite rank by suitably extending the norm from $A_0$ to $M$. If the rank is infinite, then all three notions on basis are different in general. A simple topological argument shows that the module $b(I)$ in Example 2.2 has a formal basis but does not have an orthonormal basis. Similarly, the module $c(I)$ in Example 2.3 has an orthonormal basis but does not have a formal basis. If one takes the direct sum of the modules $b(I)$ and $c(I)$, then one gets a Banach module over $A_0$ with a topological basis but without formal and orthonormal bases. By completing a Banach $A_0$-module in some way, one sees that any Banach $A_0$-module with a topological basis can be embedded as a submodule of a Banach $A_0$-module with a formal basis. Thus, the module $b(I)$ in Example 2.2 provides the largest type of Banach $A_0$-modules with a topological basis. In this paper, unless otherwise stated, we shall only consider Banach $A_0$-module of type $b(I)$. Although $c(I)$ is a submodule of $b(I)$, a continuous endomorphism of $c(I)$ cannot in general be extended to a continuous endomorphism of $b(I)$ because the topologies on $c(I)$ and $b(I)$ are different (they do have the same norm topology).

From now on, we assume that $M$ is a Banach $A_0$-module with a formal basis $\vec{e} = \{e_i | i \in I\}$ indexed by some set $I$ of at most countable cardinality. For convenience, the set $I$ will often be identified with a subset of the positive integers. In particular, we often identify our index set $I$ with the full set of positive integers whenever $M$ is of infinite rank over $A_0$. The following lemma shows which transition matrices give a new formal basis for $M$ over $A_0$, see Lemma 2.1 in [16] for a proof.

**Lemma 2.4.** *Let $M$ be a Banach $A_0$-module with a formal basis $\vec{e}$ indexed by some set $I$ of at most countable cardinality. Let $\vec{f} = \vec{e}U$ be a set of elements in $M$ indexed by $I$, where $U$ is a matrix with entries in $A_0$ whose rows and columns are both indexed by $I$. Then, $\vec{f}$ is also a formal basis of $M$ over $A_0$ if and only if $U$ is invertible over $A_0$ and the row vectors of both $U$ and $U^{-1}$ are in $c(I)$.*

We can now define the notion of $\sigma$-modules and nuclear $\sigma$-modules. Recall that $M$ has two different topologies in the infinite rank case.

**Definition 2.5.** (i). A $\sigma$-module over $A_0$ is a pair $(M, \phi)$, where $M$ is a Banach $A_0$-module with a formal basis $\vec{e} = \{e_i\}_{i \in I}$ indexed by $I$ of at most countable cardinality such that $\phi$ is a continuous $\sigma$-linear endomorphism of $M$. (ii). The $\sigma$-module $(M, \phi)$ is called nuclear if in addition, $\phi$ is also continuous when we put the norm topology on the range (this condition is called the nuclear condition). (iii). The norm of $\phi$ is defined to be $\|\phi\| = p^{-\mathrm{ord}_\pi(\phi)}$, where $\mathrm{ord}_\pi(\phi)$ is the smallest non-negative integer $k$ such that

$$\phi(M) \subset \pi^k M.$$

Since $\phi$ is continuous and $\sigma$-linear, we deduce the infinite $\sigma$-linearity of $\phi$:

$$(2.2) \qquad \phi(\sum_{i \in I} a_i e_i) = \sum_{i \in I} \sigma(a_i)\phi(e_i), \ a_i \in A_0.$$

By our definition, the nuclear condition of $\phi$ means that whenever $\{m_i\}$ is a sequence tending to zero in the topology of $M$, then we have

$$\lim_{i \to \infty} \|\phi(m_i)\| = 0.$$

From this, it is easy to check that $\phi$ is nuclear if and only if

$$(2.3) \qquad \lim_{i \to \infty} \|\phi(e_i)\| = 0.$$

The condition in (2.3) is easier to check in practice and will also be called the nuclear condition. Since the nuclear condition is a topological condition, the limit in (2.3) holds for every formal basis $\vec{e}$ of $M$.

We emphasize that our Banach $A_0$-module $M$ has a formal basis but does not have an orthonormal basis in the infinite rank case. Thus, specifying the image of $\phi$ at the basis elements $\{e_i\}$ is not enough to define a continuous $\sigma$-linear map. This is because the infinite sum $\phi(\sum a_i e_i) = \sum_i \sigma(a_i)\phi(e_i)$ may not be defined when we write out $\phi(e_i)$. In the nuclear case, our condition $\lim_i \|\phi(e_i)\| = 0$ insures that the infinite sum $\sum_i \sigma(a_i)\phi(e_i)$ is always well defined.

Intuitively, one can think of our nuclear map $\phi$ as a "family" of completely continuous linear endomorphisms of a fixed Banach space with a formal basis, parametrized by the value $x$ of $X$ at various closed points of the variety $\mathbf{A}^n$, see later definition of L-functions in this section. If one really wants to use the more familiar notion of a Banach space with an orthonormal basis, then one can think of our nuclear map $\phi$ as the continuous dual of a "family" of completely continuous linear endomorphisms of a fixed Banach space with an orthonormal basis. The word "family" is not in the usual algebraic sense, because the fibre $\phi_x$ of $\phi$ at a closed point $x$ is not a linear map but only a $\sigma$-linear map. However, a suitable power (iterate) of $\phi_x$ is linear, see the definition of L-functions in this section. In the classical case when the ground ring is $R$ or $K$, the terminology "nuclear operator" has been used (by Serre and Monsky) to denote a linear operator with certain good spectral properties. It is known that a completely continuous linear endomorphism is nuclear. Thus, we have borrowed the terminology "nuclear" for our map $\phi$ because our nuclear condition (2.3) obviously implies that $\phi$ has some good spectral properties, although we do not discuss the spectral aspect here since our map is not linear. Also, the word "nuclear" is shorter than "completely continuous".

A morphism between two $\sigma$-modules $(M, \phi)$ and $(N, \psi)$ is a continuous $A_0$-linear map of $A_0$-modules

$$\theta : M \longrightarrow N$$

such that $\theta \circ \phi = \psi \circ \theta$. In this way, the category of $\sigma$-modules is defined. In particular, it makes sense to talk about isomorphic $\sigma$-modules. It is easy to check that the usual direct sum of two $\sigma$-modules

$$(M, \phi) \oplus (N, \psi) = (M \oplus N, \phi \oplus \psi)$$

is again a $\sigma$-module. If both $\phi$ and $\psi$ are nuclear, so is the direct sum $\phi \oplus \psi$. To extend the usual tensor operations to nuclear $\sigma$-modules, we need to introduce the concepts of formal tensor product, formal symmetric powers and formal exterior powers.

**Definition 2.6.** Let $M$ and $N$ be two Banach $A_0$-modules. Their formal tensor product, denoted by $M \bar{\otimes} N$, is the completion of the usual tensor product with respect to the topology induced by requiring the usual tensor product $P \otimes Q$, where $P$ is an open neighborhood of the origin in $M$ and $Q$ is an open neighborhood of the origin in $N$, to be a basis of open neighborhoods of the origin. The norm on the formal tensor product $M \otimes N$ is defined as follows. If $c$ is an element of the usual algebraic tensor product, then

$$\|c\| = \min_{c = \Sigma_j m_j \otimes n_j} \max_j \{\|m_j\| \cdot \|n_j\|\},$$

where the sum runs over all ways that $c$ can be written as a finite sum $c = \sum_j m_j \otimes n_j$. For a general element $c$ in the formal tensor product $M \otimes N$, we define

$$\|c\| = \min_{\{c_i\}} \liminf_i \|c_i\|,$$

where $\{c_i\}$ runs over all sequences approaching to $c$ such that each $c_i$ is in the usual algebraic tensor product of $M$ and $N$.

Note that we are using the same tensor notation to denote the formal tensor product. This should not cause confusion, since all tensor products in this paper are formal tensor products. If $\vec{e} = \{e_i\}_{i \in I}$ is a formal basis of $M$ and $\vec{f} = \{f_j\}_{j \in J}$ is a formal basis of $N$, then one checks that $\vec{e} \otimes \vec{f} = \{e_i \otimes f_j : i \in I, j \in J\}$ is a formal basis of $M \bar{\otimes} N$. Formal symmetric products and formal exterior products can be defined and dealt with similarly.

**Definition 2.7.** Let $(M, \phi)$ be a $\sigma$-module with a formal basis $\vec{e}$. Let $(N, \psi)$ be a $\sigma$-module with a formal basis $\vec{f}$. Let

$$\phi(e_i) = \sum_k a_{ki} e_i, \ \psi(f_j) = \sum_k b_{kj} f_j.$$

The formal tensor product of $(M, \phi)$ and $(N, \psi)$ is the $\sigma$-module $(M \bar{\otimes} N, \phi \otimes \psi)$, where $M \otimes N$ is the formal tensor product of the two Banach $A_0$-modules $M$ and

$N$, the map $\phi \otimes \psi$ is determined by the relation

$$\begin{aligned}
\phi \otimes \psi(e_i \otimes f_j) &= \phi(e_i) \otimes \psi(f_j) \\
&= (\sum_{k_1} a_{k_1 i} e_{k_1}) \otimes (\sum_{k_2 j} b_{k_2 j} f_{k_2}) \\
&= \sum_{k_1, k_2} a_{k_1 i} b_{k_2 j}(e_{k_1} \otimes f_{k_2}).
\end{aligned}$$

If both $\phi$ and $\psi$ are nuclear, it is easy to check that

$$\lim_{i+j \to \infty} \|\phi \otimes \psi(e_i \otimes f_j)\| = 0.$$

Thus, the formal tensor product $(M \otimes N, \phi \otimes \psi)$ is also a nuclear $\sigma$-module.

Let $(M, \phi)$ be a $\sigma$-module with a formal basis $\vec{e}$. In a similar way, one can define formal exterior powers and formal symmetric powers. Let $1 \leq i < \infty$. The $i$-th formal exterior power $(\wedge^i M, \wedge^i \phi)$ is the $\sigma$-module with the formal basis

$$\wedge^i \vec{e} = \{\cdots, e_{k_1} \wedge \cdots \wedge e_{k_i}, \cdots\}, \ k_1 < k_2 < \cdots < k_i$$

and with the continuous $\sigma$-linear endomorphism $\wedge^i \phi$ defined by

$$\wedge^i(\phi)(e_{k_1} \wedge \cdots \wedge e_{k_i}) = \phi(e_{k_1}) \wedge \cdots \wedge \phi(e_{k_i}).$$

For $i = 0$, we define $(\wedge^0 M, \wedge^0 \phi)$ to be the trivial rank one $\sigma$-module $(A_0, \sigma)$. Similarly, the $i$-th formal symmetric power $(\mathrm{Sym}^i M, \mathrm{Sym}^i \phi)$ is the $\sigma$-module with the formal basis

$$\mathrm{Sym}^i \vec{e} = \{\cdots, e_{k_1} e_{k_2} \cdots e_{k_i}, \cdots\}, \ k_1 \leq k_2 \leq \cdots \leq k_i$$

and with the continuous $\sigma$-linear endomorphism $\mathrm{Sym}^i \phi$ defined by

$$\mathrm{Sym}^i(\phi)(e_{k_1} \cdots e_{k_i}) = \phi(e_{k_1}) \cdots \phi(e_{k_i}).$$

We shall use the convention that $\mathrm{Sym}^k \phi = 0$ for all negative integers $k < 0$. If $\phi$ is nuclear, one checks that its formal symmetric powers and formal exterior powers are also nuclear. We obtain

**Lemma 2.8.** *The category of $\sigma$-modules over $A_0$ (resp. nuclear $\sigma$-modules over $A_0$) is closed under direct sums, formal tensor product, formal symmetric powers and formal exterior powers.*

Another useful construction is the formal symmetric algebra. Let $M$ be a Banach $A_0$-module. We define the **formal symmetric algebra** of $M$, denoted by $A_0[[M]]$, to be the completion of the ring $\oplus_{i \geq 0} \mathrm{Sym}^i M$ with respect to the ideal $\oplus_{i \geq 1} \mathrm{Sym}^i M$, where $\mathrm{Sym}^i M$ is the $i$-th formal symmetric power of $M$. An element $m \in A_0[[M]]$ can be written uniquely as a series

$$m = \sum_{i \geq 0} m_i, \ m_i \in \mathrm{Sym}^i M.$$

The ring $A_0[[M]]$ is also a Banach $A_0$-module with respect to the norm defined by

$$\|m\| = \max_i \|m_i\|.$$

If $\vec{e} = \{e_1, e_2, \cdots\}$ is a formal basis of $M$, then

$$\{e_{\ell_1} e_{\ell_2} \cdots e_{\ell_i}; i \geq 0\}$$

is a formal basis for $A_0[[M]]$. If $\phi$ is a continuous $\sigma$-linear map on $M$, then $\phi$ induces a continuous $\sigma$-linear map $\mathrm{Sym}(\phi)$ on $A_0[[M]]$. It should be remarked that if $\phi$ is nuclear on $M$, then $\mathrm{Sym}(\phi)$ is in general **not** nuclear. The easiest counter-example is to take the identity $\sigma$-module $(M, \phi) = (A_0, \sigma)$.

The continuous $\sigma$-linear map $\phi$ in a $\sigma$-module $(M, \phi)$ can also be viewed as a continuous $A_0$-linear map $\phi : M^\sigma \to M$ of Banach $A_0$-modules, where $M^\sigma$ is the continuous pull back of $M$ by $\sigma : A_0 \to A_0$. That is,

$$M^\sigma = \{\sum_i b_i \otimes e_i \mid b_i \in A_0\}$$

with the formal basis $\{1 \otimes e_i\}_{i \in I}$. For $a_i \in A_0$, we have the relation

$$\sum \sigma(a_i) \otimes e_i = \sum_i a_i e_i.$$

The map $\phi$ is called the Frobenius map of the $\sigma$-module $M$. We shall sometimes just write $M$ or $\phi$ for the pair $(M, \phi)$. A $\sigma$-module $\phi$ is called a **unit root** $\sigma$-module if the $A_0$-linear map $\phi : M^\sigma \to M$ is an isomorphism of $A_0$-modules. The nuclear condition (2.3) implies that any unit root nuclear $\sigma$-module $\phi$ must be of finite rank. If $M$ is of finite rank, the two topologies on $M$ agrees and hence the nuclear condition (2.3) is automatically satisfied. Thus, in finite rank case, we may drop the word "nuclear" and simply say a $\sigma$-module.

Given a $\sigma$-module $(M, \phi)$, we define its dual $A_0$-module by

$$M^\vee = \mathrm{Hom}^{\mathrm{cont}}_{A_0}(M, A_0).$$

This is the set of continuous $A_0$-linear maps $f$ from $M$ to $A_0$, i.e., $A_0$-linear maps satisfying

$$f(\sum_i a_i e_i) = \sum_i a_i f(e_i), \ a_i \in A_0.$$

A given $A_0$-linear map $f$ from $M$ to $A_0$ is continuous if and only if $f$ has the property

$$\lim_i f(e_i) = 0.$$

The module $M^\vee$ is a Banach $A_0$-module if we define the norm on it by

$$\|f\| = \max_{m \in M} \|f(m)\|.$$

Let $\vec{e}^\vee = \{e_1^\vee, e_2^\vee, \cdots\}$ be the dual basis

$$e_j^\vee(e_i) = \begin{cases} 1, & \text{if } j = i, \\ 0. & \text{if } j \neq i, \end{cases}$$

Then, one checks that

$$M^\vee = \{\sum_i a_i e_i^\vee \mid a_i \in A_0, \ \lim_i \|a_i\| = 0\}$$

is the Banach $A_0$-module with orthonormal basis $\vec{e}^\vee$. The dual map

$$\phi^\vee : M^\vee \longrightarrow M^{\sigma\vee}, \ \phi^\vee(f) = f \circ \phi, \ f \in M^\vee$$

is then continuous and $A_0$-linear, where

$$M^{\sigma\vee} = \{\sum a_i \otimes e_i^\vee \mid a_i \in A_0, \ \lim_i \|a_i\| = 0, \ \sigma(a_i) \otimes e_i^\vee = a_i e_i^\vee\}.$$

If $\phi$ is nuclear, then as an $A_0$-linear map between two different spaces, $\phi^\vee$ is completely continuous as defined in Coleman [2] where the case of completely continuous linear endomorphisms (between the same space, not two different spaces) was studied in detail. The above dual $(M^\vee, \phi^\vee)$ is not a $\sigma$-module yet by our definition since the arrow direction for the map $\phi^\vee$ is not the right one and furthermore $M^\vee$ does not have a formal basis. In the case that $\phi$ is invertible, the contragredient (the inverse of $\phi^\vee$) of $\phi$ can be used to define a $\sigma$-module on the completion of $M^\vee$. We shall not need this construction. All we need in this paper is the completely continuous property of $\phi^\vee$.

We now turn to discussing some "convergence" conditions of $\phi$ as a "power series" in terms of $X$, generalizing the convergent and overconvergent conditions for an element in $A_0$. Finer explicit description about the nuclear condition of $\phi$ will be discussed in section 6 in connection with the Hodge-Newton decomposition. Let $B(X)$ be the matrix of a $\sigma$-module $(M, \phi)$ under a formal basis $\vec{e}$. We can always write

$$B(X) = \sum_{u \in \mathbf{Z}_{\geq 0}^n} B_u X^u,$$

where the coefficients $B_u$ are matrices with entries in $R$.

**Definition 2.9.** We say that $(M, \phi)$ is **convergent** if there is a formal basis $\vec{e}$ such that the matrix $B(X)$ with respect to $\vec{e}$ is convergent, that is, $B(X)$ satisfies the condition

$$(2.4) \qquad \lim_{|u| \to \infty} \operatorname{ord}_\pi B_u = \infty.$$

Similarly, we say that $(M, \phi)$ is **overconvergent** if there is a formal basis $\vec{e}$ such that the matrix $B(X)$ with respect to $\vec{e}$ is overconvergent. Namely, $B(X)$ satisfies the condition

$$(2.5) \qquad \lim_{|u| \to \infty} \inf \frac{\operatorname{ord}_\pi B_u}{|u|} > 0.$$

Note that our definition of the above convergence properties depends on the choice of the formal basis $\vec{e}$. It may happen that the matrix $B(X)$ under one formal basis is overconvergent but the matrix under a new formal basis is NOT overconvergent. Our definition is that as long as there is one formal basis for which the matrix is overconvergent, then the $\sigma$-module $(M, \phi)$ is called overconvergent. When we say a formal basis $\vec{e}$ of a convergent (resp. overconvergent) $\phi$, we always assume that the matrix of $\phi$ with respect to $\vec{e}$ is convergent (resp. overconvergent). Alternatively, we may say that the basis $\vec{e}$ is convergent (resp. overconvergent) for $(M, \phi)$ if the matrix of $\phi$ with respect to $\vec{e}$ is convergent (resp. overconvergent). In the case of finite rank, every $\sigma$-module is convergent; an overconvergent $\sigma$-module is a $\sigma$-module which can be defined over $A$ and then by base extension $A \to A_0$. In the case of infinite rank, not every $\sigma$-module is convergent. It may happen that the matrix $B(X) = (b_{i,j}(X))$ has all of its entries $b_{i,j}(X)$ in $A_0$ without satisfying the uniformity requirement (2.4). Similar remarks apply to the overconvergent situation. It may happen that the matrix $B(X) = (b_{i,j}(X))$ has all of its entries $b_{i,j}(X)$ in $A$ without satisfying the uniformity requirement (2.5).

The matrix $B(X)$ of a $\sigma$-module with respect to a formal basis can be further written in the form

$$B(X) = \sum_{u \in \mathbf{Z}_{\geq 0}^n} B_u X^u, \quad B_u = (b_{w_1, w_2}(u)),$$

where $w_1$ denotes the row index and $w_2$ denotes the column index of the constant matrix $B_u$. Conversely, such a matrix power series $B(X)$ is the matrix of some overconvergent nuclear $\sigma$-module under some formal basis if and only if $B(X)$ satisfies the following two conditions. First, we have the overconvergent condition for $B(X)$ as a power series in $X$:

$$\lim_{|u| \to \infty} \inf \frac{\mathrm{ord}_\pi B_u}{|u|} > 0.$$

Second, we have the nuclear condition for $B(X)$ to be the matrix of $\phi$ under a formal basis. This nuclear condition can be rephrased as follow: For any positive number $C > 0$, there is an integer $N_C > 0$ such that for all column indices $w_2 > N_C$, we have

$$\mathrm{ord}_\pi b_{w_1, w_2}(u) \geq C,$$

uniformly for all $u$ and all $w_1$.

As stated in Lemma 2.8, the category of $\sigma$-modules (resp. nuclear $\sigma$-modules) is closed under direct sum, formal tensor product, formal symmetric product and formal exterior product. The same statement holds for the subcategory of convergent (resp. overconvergent) nuclear $\sigma$-modules. The category of convergent nuclear $\sigma$-modules includes all $\sigma$-modules of finite rank. In the case that the finite rank map $\phi : M^\sigma \to M$ is an isogeny, namely when $\phi : M^\sigma \otimes \mathbf{Q} \to M \otimes \mathbf{Q}$ is an isomorphism of finite rank $A_0 \otimes \mathbf{Q}$-modules, the $\sigma$-module $(M, \phi)$ more or less corresponds to an F-crystal on the affine $n$-space $\mathbf{A}^n$. If, in addition, the finite rank map $\phi : M^\sigma \to M$ is an isomorphism of $A_0$-modules, the $\sigma$-module $(M, \phi)$ corresponds exactly to a unit root F-crystal on $\mathbf{A}^n$. However, the category of finite rank overconvergent unit root $\sigma$-modules is much larger than the category of finite rank overconvergent unit root F-crystals. The notion of overconvergent nuclear $\sigma$-modules of finite rank arises naturally from various relative $p$-adic (but not $p$-adic étale) cohomology groups of a family of varieties over finite fields. The notion of overconvergent nuclear $\sigma$-modules of infinite rank already arises if one works on the chain level of $p$-adic cohomology groups of a family of varieties over a finite field.

From now on, in order to define L-function, we shall only work with a **nuclear $\sigma$-module** $(M, \phi)$ and sometimes identify it with its matrix $B(X)$ under a formal basis. One can attach an L-function to any nuclear $\sigma$-module $\phi$. For this purpose, we need to recall the notion of the Techmüller lifting of a geometric point $x \in \mathbf{A}^n(\bar{\mathbf{F}}_q)$ defined with respect to $\sigma$. Let $\bar{x} \in \mathbf{A}^n$ be a geometric point of degree $d$ over $\mathbf{F}_q$. That is, the point $\bar{x}$ is a surjective $\mathbf{F}_q$-algebra homomorphism in $\mathrm{Hom}_{\mathbf{F}_q}(\mathbf{F}_q[X], \mathbf{F}_{q^d})$. By a theorem of Monsky-Tate, there is a unique surjective $R$-algebra homomorphism $x \in \mathrm{Hom}_R(A_0, R_d)$ whose reduction mod $\pi$ is $\bar{x}$ such that $x \circ \sigma^d = x$, where $R_d$ denotes the unramified extension of $R$ with residue field $\mathbf{F}_{q^d}$. The element $x$ is called the Teichmüller lifting of $\bar{x}$ defined with respect to $\sigma$. The map $\sigma$ on $A_0$ specializes to an action on $R_d$ via the element $x \in \mathrm{Hom}_R(A_0, R_d)$. This action on $R_d$ coincides with the standard Frobenius action on $R_d$.

The fibre $M_x$ of $M$ at $x$ is a Banach module with a formal basis over $R[x]$, where $R[x] \ (= R_d)$ is the ring obtained from $R$ by adjoining the coordinates $(x_1, \cdots, x_n)$ of the point $x$ defined by $x_i = x(X_i)$. That is, via the map $x$, $R[x]$ becomes a Banach $A_0$-algebra and

$$M_x = M \otimes R[x]$$

is the formal tensor product of $M$ and $R[x]$ via $x$. It has a formal basis as a Banach $R[x]$-module but has no topological basis as a Banach $A_0$-module. The fibre map

$$\phi_x : M_x \longrightarrow M_x$$

is $\sigma$-linear (not linear in general). But its $d$-th iterate $\phi_x^d$ is $R[x]$-linear (since $\sigma^d(x) = x$) and can be defined over $R$. Furthermore, $\phi_x^d$ is the dual of a completely continuous operator in the sense of Serre [11]. This implies that the characteristic series $\det(I - \phi_x^d T)$ (or the Fredholm determinant) of $\phi_x^d$ acting on $M_x$ (which can be descended from $R[x]$ to $R$) is well defined with coefficients in $R$. It is a $p$-adic entire function and is independent of the choice of the geometric point $\bar{x}$ associated to the closed point over $\mathbf{F}_q$ containing $\bar{x}$. Thus, it makes sense to talk about the characteristic power series of $\phi$ at a closed point.

**Definition 2.10.** Let $(M, \phi)$ be a nuclear $\sigma$-module over $A_0$. The L-function of $(M, \phi)$ is defined by the Euler product

$$L(\phi, T) = \prod_{\bar{x} \in \mathbf{A}^n/\mathbf{F}_q} \frac{1}{\det(I - \phi_x^{\deg(\bar{x})} T^{\deg(\bar{x})})},$$

where $\bar{x}$ runs over the closed points of $\mathbf{A}^n/\mathbf{F}_q$ and $x$ is the Teichmüller lifting of $\bar{x}$ defined with respect to $\sigma$.

Note that each of our Euler factors is the reciprocal of an entire power series, not necessarily a polynomial. If $B(X)$ is the matrix of $\phi$ with respect to a formal basis, then one checks that

$$\begin{aligned}
L(\phi, T) &= L(B(X), T) \\
&= \prod_{\bar{x} \in \mathbf{A}^n/\mathbf{F}_q} \frac{1}{\det(I - B(x)B(x^\sigma) \cdots B(x^{\sigma^{\deg(\bar{x})-1}}) T^{\deg(\bar{x})})},
\end{aligned}$$

where the notation $x^{\sigma^i}$ denotes $\sigma^i(x)$. This shows that the L-function is a power series with coefficients in $R$. Thus, it is trivially meromorphic in the open unit disk $|T|_\pi < 1$. If $\phi$ is convergent (which is always the case if $M$ is of finite rank), it can be shown that the L-function is meromorphic on the closed unit disk $|T|_\pi \leq 1$. The generalized Monsky trace formula shows that the L-function $L(\phi, T)$ is $p$-adic meromorphic if $\phi$ is overconvergent and of finite rank. Our first task in this paper is then to further extend the Monsky trace formula from finite rank case to infinite rank case. This will be done in section 5. In particular, it implies

**Theorem 2.11.** *Let $(M, \phi)$ be an overconvergent nuclear $\sigma$-module. Then, the L-function $L(\phi, T)$ is $p$-adic meromorphic.*

Actually, what we will need in this paper is a stronger family version of Theorem 2.11 as given later in section 5. In fact, we will go further and work out the effective estimates involved.

In order to fully describe Dwork's conjecture in our language, we need to define a slightly more general L-function, which could be called higher power L-function. If $\phi$ is a nuclear $\sigma$-module and $k$ is a positive integer, we define the $k$-th power $\phi^k$ to be $\phi$ iterated $k$-times. The pair $(M, \phi^k)$ can be viewed as a $\sigma^k$-module but it is **not** a $\sigma$-module as $\phi^k$ is $\sigma^k$-linear not $\sigma$-linear in general. The L-function of the

$k$-th power $\phi^k$ is defined to be

$$L(\phi^k, T) = \prod_{\bar{x} \in \mathbf{A}^n/\mathbf{F}_q} \frac{1}{\det(I - \phi_x^{k \deg(\bar{x})} T^{\deg(\bar{x})})}.$$

If $B(X)$ is the matrix of $\phi$ with respect to a formal basis, then one checks that

$$
\begin{aligned}
L(\phi^k, T) &= \prod_{\bar{x} \in \mathbf{A}^n/\mathbf{F}_q} \frac{1}{\det(I - B(x)B(x^\sigma) \cdots B(x^{\sigma^{k \deg(\bar{x})-1}}) T^{\deg(\bar{x})})} \\
&= \prod_{\bar{x} \in \mathbf{A}^n/\mathbf{F}_q} \frac{1}{\det(I - (B(x)B(x^\sigma) \cdots B(x^{\sigma^{\deg(\bar{x})-1}}))^k T^{\deg(\bar{x})})}.
\end{aligned}
$$

In other words, $L(\phi^k, T)$ is obtained from $L(\phi, T)$ by raising each reciprocal root $\alpha$ of the Euler factors to its $k$-th power $\alpha^k$.

Let $(M_1, \phi_1)$ and $(M_2, \phi_2)$ be two nuclear $\sigma$-modules. Let $k_1$ and $k_2$ be two positive integers. Although their powers $\phi_1^{k_1}$ and $\phi_2^{k_2}$ do not have a $\sigma$-module structure in general, we can define the L-function of the "tensor product" of these powers as follows:

$$L(\phi_1^{k_1} \otimes \phi_2^{k_2}, T) = \prod_{\bar{x} \in \mathbf{A}^n/\mathbf{F}_q} \frac{1}{\det(I - \phi_{1x}^{k_1 \deg(\bar{x})} \otimes \phi_{2x}^{k_2 \deg(\bar{x})} T^{\deg(\bar{x})})}.$$

Namely, the Euler factor at $x$ of $L(\phi_1^{k_1} \otimes \phi_2^{k_2}, T)$ is the "tensor product" of the corrosponding Euler factors at $x$ of $L(\phi_1^{k_1}, T)$ and $L(\phi_2^{k_2}, T)$. The above definition works for negative integer $k$ as well provided that $\phi$ is injective fibre by fibre and of finite rank. If $(M, \phi)$ is of rank one, then $L(\phi^k, T)$ is just the L-function $L(\phi^{\otimes k}, T)$ of the $k$-th tensor (or symmetric) power of $(M, \phi)$. If $(M, \phi)$ is of rank greater than one, then the $k$-th power $\phi^k$ $(k > 1)$ does not have a $\sigma$-module interpretation as matrix multiplication is non-commutative in general. Nevertheless, $\phi^k$ can be viewed as a virtual $\sigma$-module, namely, an integral combination of $\sigma$-modules with both positive and negative coefficients. In fact, we have the following simple but basic decomposition formula for L-functions, see equality (4.4) in [16].

**Lemma 2.12.** *Let $\phi$ be a nuclear $\sigma$-module. Let $k$ be a positive integer. Then,*

$$L(\phi^k, T) = \prod_{i \geq 1} L(\mathrm{Sym}^{k-i} \phi \otimes \wedge^i \phi, T)^{(-1)^{i-1} i}.$$

Since $\mathrm{Sym}^{k-i} \phi = 0$ for $i > k$, Lemma 2.12 implies that $L(\phi^k, T)$ can be written as an finite alternating product of L-functions of various $\sigma$-modules. It is known that the L-function $L(\phi^k, T)$ is not meromorphic in general if $\phi$ is not overconvergent. Dwork's conjecture concerns the meromorphic continuation of $L(\phi_{\mathrm{unit}}^k, T)$ for certain unit root $\sigma$-modules $\phi_{\mathrm{unit}}$ which are not overconvergent but arise in a natural way from the Hodge-Newton decomposition of an overconvergent $\sigma$-module $\phi$. Since $\phi_{\mathrm{unit}}$ is not overconvergent, meromorphic continuation in this case is much deeper than Theorem 2.11. To describe and prove the full strength of our results, we need more preparations.

## 3. Uniform and strong families

Since it is also our objective to understand how the arithmetic family of L-functions $L(\phi^k, T)$ varies when the integer parameter $k$ varies, we first review some notions from [16] about the variation of meromorphic functions.

Let $S$ be a subset of a complete topological space. Let $\bar{S}$ denote the closure of $S$. Namely, $\bar{S}$ is the union of $S$ and its limit points. In our applications, the space $S$ will be an infinite subset of the integers $\mathbf{Z}$, sometimes with the induced $p$-adic topology, and sometimes with the sequence topology with $\infty$ as its unique limit point not in it. Recall that the sequence topology of an infinite subset $S \subset \mathbf{Z}$ means that the open sets of $S$ are the finite subsets of $S$ as well as compliments of finite subsets of $S$. One checks that for the sequence topology, a Cauchy sequence of $S$ is the same as a continuous sequence of $S$. The set $S$ will be our parameter space for which our parameter $k$ varies. We shall consider a family of power series in $R[[T]]$ of the following form:

$$f(k, T) = \sum_{m \geq 0} f_m(k) T^m, \ f_0(k) = 1, \ f_m(k) \in R,$$

where each coefficient $f_m(k)$ is a function from $S$ to $R$. Note that the constant term is always 1. Furthermore, the ring $R$ is always complete and compact with respect to the $p$-adic topology of $R$. The power series ring $R[[T]]$ is a Banach $R$-algebra under the Gauss norm

$$\| \sum a_m T^m \| = \max_m \| a_m \|.$$

It has the formal basis $\{T^i : i \geq 0\}$. The family $f(k, T)$ is called **uniformly continuous** in $k$ if the norm $\| f(k_1, T) - f(k_2, T) \|$ is uniformly small whenever $k_1$ is close to $k_2$. The family of functions $f(k, T)$ is called **uniformly entire** if

$$(3.1) \qquad \lim_{m \to \infty} \inf \frac{\inf_k \operatorname{ord}_\pi f_m(k)}{m} = +\infty.$$

The family $f(k, T)$ is called **uniformly meromorphic** if $f(k, T)$ can be written as a quotient

$$f(k, T) = \frac{f_1(k, T)}{f_2(k, T)}$$

of two families, where both $f_1(k, T)$ and $f_2(k, T)$ are uniformly entire. It is clear that the product and quotient (if the denominator is not the zero family) of two uniformly meromorphic families (parameterized by the same parameter $k$) are still uniformly meromorphic.

If the set $S$ consists of a sequence of elements $\{k_i\}$, we shall call the family $f(k_i, T)$ a sequence of functions. The following result is an immediate consequence of the above defining inequality in (3.1).

**Lemma 3.1.** *Let $f(k_i, T)$ be a uniformly continuous sequence of power series. Assume that the sequence $f(k_i, T)$ is uniformly entire. If $k = \lim_i k_i$ exists in $\bar{S}$, then the limit function*

$$f(k, T) = \lim_i f(k_i, T)$$

*exists and it is also entire.*

**Corollary 3.2.** *If $f(k, T)$ is a uniformly continuous family of uniformly entire functions for $k \in S$, then $f(k, T)$ extends uniquely to a uniformly continuous family of uniformly entire functions for $k \in \bar{S}$ (the closure of $S$).*

To obtain similar result for meromorphic families, we need to introduce another notion about families of functions.

**Definition 3.3.** A family $f(k, T)$ is called a **strong family** if $f(k, T)$ is the quotient of two uniformly continuous families of uniformly entire functions.

There is probably a better choice of terminology for the notion of a strong family. We have tried to avoid the weaker notion of a uniformly continuous family of uniformly meromorphic functions, which simply means a uniformly continuous family which is also uniformly meromorphic. In such a family, it is not immediately obvious (although probably true) that the limit function is still a meromorphic function. In our stronger definition, this property is built in because of Corollary 3.2. Thus, we have

**Corollary 3.4.** *If $f(k,T)$ is a strong family for $k \in S$, then $f(k,T)$ extends uniquely to a strong family for $k \in \bar{S}$ (the closure of $S$).*

It is clear that the product and quotient (if denominator is not the zero family) of two strong families for $k \in S$ are again a strong family.

Our main concern will be the family of L-functions arising from a family of nuclear $\sigma$-modules. We want to know when such a family of L-functions are uniformly meromorphic and when it is a strong family. Before doing so, we need to define the notion of a uniform family of $\sigma$-modules. Let $B(k,X)$ be the matrix of a family of nuclear $\sigma$-modules $(M_k, \phi_k)$ parametrized by a parameter $k$. There may be no relations among the $\phi_k$ for different $k$. In particular, the rank of $(M_k, \phi_k)$ could be totally different (some finite and other infinite, for instance) as $k$ varies. Write

$$B(k,X) = \sum_{u \in \mathbf{Z}_{\geq 0}^n} B_u(k) X^u,$$

where each coefficient matrix $B_u(k)$ has entries in $R$:

$$B_u(k) = (b_{w_1, w_2}(u,k)), \quad b_{w_1, w_2}(u,k) \in R,$$

where $w_1$ is the row index of $B_u(k)$ and $w_2$ is the column index of $B_u(k)$.

**Definition 3.5.** A family $(M_k, \phi_k)$ is called a family of **uniformly** overconvergent nuclear $\sigma$-modules if it has a family of matrices $B(k,X)$ satisfying the following two conditions. First, we have the uniform overconvergent condition:

$$\lim_{|u| \to \infty} \inf \frac{\inf_k \operatorname{ord}_\pi B_u(k)}{|u|} > 0.$$

Second, we have the uniform nuclear condition. That is, for any positive number $C > 0$, there is an integer $N_C > 0$ such that for all column number $w_2 > N_C$, we have

$$\operatorname{ord}_\pi b_{w_1, w_2}(u,k) \geq C,$$

uniformly for all $u, k, w_1$.

Note that the uniform nuclear condition is automatically satisfied if the rank of $(M_k, \phi_k)$ is uniformly bounded. This is the case if each $\sigma$-module $M_k$ has the same finite rank for every $k$. On the other hand, the uniform overconvergent condition is not automatic even if $M_k$ has the same finite rank for every $k$.

**Definition 3.6.** A family $(M_k, \phi_k)$ of nuclear $\sigma$-modules with the same underlying module $M_k = M$ for $k \in S$ is called uniformly continuous in $k$ if $||\phi_{k_1} - \phi_{k_2}||$ is uniformly small whenever $k_1$ is uniformly close to $k_2$.

In applications, the underlying module $M_k$ for each $k$ may be different. In this case, it is often useful to extend each nuclear $\sigma$-module $(M_k, \phi_k)$ to a single underlying module $M$ independent of $k$ so that we have

$$M_k \oplus N_k = M.$$

The map $\phi_k$ extends to $M$ by requiring $\phi_k(\alpha) = 0$ for $\alpha \in N_k$. This will not change the L-function. That is, $L(\phi_k, T)$ is the same whether we view $\phi_k$ as acting on the bigger module $M$ or on the smaller submodule $M_k$. In this way, we can talk about the uniform continuity as defined in Definition 3.6. The following uniform result is of basic importance to our investigation.

**Theorem 3.7.** *Let $\phi_k$ be a family of nuclear $\sigma$-modules parameterized by $k$. Assume that the family $\phi_k$ is uniformly overconvergent. Then, the family of L-functions $L(\phi_k, T)$ is uniformly meromorphic. If, in addition, the family $\phi_k$ is uniformly continuous, then the family $L(\phi_k, T)$ of L-functions parametrized by $k$ is a strong family.*

This result is an extension of Theorem 2.11 from a single nuclear overconvergent $\sigma$-module to a uniform family of nuclear overconvergent $\sigma$-modules. In the special case when $\sigma(x) = x^q$, Theorem 3.7 was already proved using the extended Dwork trace formula, see Theorem 5.7 in [16]. We shall prove Theorem 3.7 in the general case in next two sections, after we establish the extended Monsky trace formula.

## 4. Explicit estimates on Dwork operators

Recall that $A$ is the overconvergent subring of $A_0$ and $\sigma$ is a fixed $R$-algebra endomorphism of $A_0$ lifting the $q$-th power Frobenius map, where $R$ is a complete discrete valuation ring of characteristic zero with uniformizer $\pi$ and residue field $\mathbf{F}_q$. As we assumed throughout, the map $\sigma$ is defined over $A$. That is, $\sigma(X_i)$ is overconvergent. It is clear that $\sigma(A_0)$ is a subring of $A_0$ and $\sigma$ is an injective map. For a lattice point $u \in \mathbf{Z}_{\geq 0}^n$, we use $u_i$ to denote its $i$-th coordinate. That is, $u = (u_1, \cdots, u_n)$. In this section, we make explicit the various estimates involving a Dwork operator.

**Definition 4.1.** A Dwork operator $\Theta$ on $A_0$ is an $R$-linear endomorphism of the ring $A_0$ which is $\sigma^{-1}$-linear in the sense that
$$\Theta(\sigma(a_1)a_2) = a_1\Theta(a_2), \quad a_1, a_2 \in A_0.$$
The operator $\Theta$ is called overconvergent if $\Theta$ is stable on the overconvergent subring $A$, that is, if
$$\Theta(A) \subset A.$$

Our main concern will be the class of overconvergent Dwork operators. To construct some examples of Dwork operators, we need the following result whose proof is straightforward.

**Lemma 4.2.** *The ring $A_0$ is a free $\sigma(A_0)$-module of rank $q^n$ with a set of generators $\{X^u\}$, where $0 \leq u_i < q$.*

For any given $f \in A_0$, Lemma 4.2 shows that there are well defined elements $\Theta_u(f)$ in $A_0$ for $0 \leq u_i < q$ such that
$$f = \sum_{0 \leq u_i < q} \Theta_u(f)^\sigma X^u.$$
For each $u$ with $0 \leq u_i < q$, the map $\Theta_u(f)$ then defines a Dwork operator on $A_0$. These $\Theta_u$ are called the basic Dwork operators of $A_0$. We will show that they are in fact overconvergent. That is, $\Theta_u$ takes $A$ to $A$. One can get more Dwork operators by composing $\Theta_u$ with an $A_0$-linear endomorphism of $A_0$.

**Definition 4.3.** For rational numbers $b > 0$ and $c$, we define a subspace of $A_0$ by

$$L(b, c) = \{ \sum_{v \in \mathbf{Z}_{\geq 0}^n} a_v X^v | a_v \in R, \ \mathrm{ord}_\pi(a_v) \geq b|v| + c\}.$$

The space $L(b, c)$ is simply the $R$-submodule of $A_0$ consisting of those power series $f(X)$ which converge on the polydisk $\mathrm{ord}_\pi(X_i) > -b$ and satisfy $\mathrm{ord}_\pi(f(X)) \geq c$ there. Since $b > 0$, the space $L(b, c)$ is an $R$-submodule of the overconvergent subring $A$. Although $A$ itself is not complete, each $L(b, c)$ is easily seen to be complete with respect to the $\pi$-adic topology. The space $L(b, c)$ is closed in $A_0$. Furthermore, we have the following simple properties:

$$A = \bigcup_{b > 0, c} L(b, c),$$

$$L(b, c_1) \cdot L(b, c_2) \subset L(b, c_1 + c_2).$$

In particular, $L(b, 0)$ is a ring and each $L(b, c)$ is an $L(b, 0)$-module. The space

$$L(b, 0) \otimes K = L(b, c) \otimes K$$

is a $p$-adic Banach space over the field $K$ with the formal basis $\{\pi^{[b|u|]} X^u\}$. The above definition works for real numbers $b > 0$ and $c$ as well, not just for rational numbers. We choose rational numbers because by going to a totally ramified finite extension of $R$ if necessary, we may assume that $b$ is a positive integer. This will allow us to replace the integer part $[b|u|]$ by $b|u|$ and to use the simpler looking orthonormal basis $\{\pi^{b|u|} X^u\}$.

We shall show that any overconvergent Dwork operator $\Theta$ acting on $A_0$ induces a completely continuous operator on $L(b, c) \otimes K$ for sufficiently small $b$. This implies that the Fredholm determinant

$$\det(I - T\Theta|A_0) = \det(I - T\Theta|L(b, c) \otimes K)$$

is well defined and is $p$-adically entire. It is independent of the choice of $b$ and $c$.

Since $\sigma(X_i) \in A$, we can write

$$\sigma(X_i) = X_i^q + \pi f_i, \ f_i \in A.$$

Choose and fix a positive rational number $b_\sigma > 0$ depending on $\sigma$ such that

$$\pi f_i \in L(b_\sigma, 0)$$

for all $1 \leq i \leq n$. If $b$ is a rational number with $0 < b \leq b_\sigma$, then

$$\pi f_i \in L(b_\sigma, 0) \subset L(b, 0) \subset L(b, -qb)$$

and trivially

$$X^q \in L(b, -qb).$$

This proves the following

**Lemma 4.4.** *For all rational numbers $0 < b \leq b_\sigma$ and all $1 \leq i \leq n$, we have*

$$\sigma(X_i) \in L(b, -qb).$$

**Theorem 4.5.** *Let $b$ and $c$ be rational numbers with $0 < b \leq b_\sigma$. For each $u$ with $0 \leq u_i < q$, we have*

$$\Theta_u(L(b, c)) \subseteq L(qb, c).$$

*In particular, each $\Theta_u$ is overconvergent.*

*Proof.* First, we assume that $f$ is the monomial $X^v$. We can write $v = qs + u$ with $0 \leq u_i < q$. Then

$$X^v = \sigma(X^s)X^u + X^u(X^{qs} - \sigma(X^s)).$$

Since

$$\sigma(X_i) \subset L(b, -qb),$$

we deduce

$$\sigma(X^s) \subset L(b, -qb|s|).$$

On the other hand,

$$\sigma(X^s) \equiv X^{qs} \ (\text{mod } \pi), \quad X^{qs} \in L(b, -qb|s|).$$

It follows that

$$\begin{aligned}
X^{qs} - \sigma(X^s) &\in \pi A \cap L(b, -qb|s|) \\
&\subseteq \pi L(b, -qb|s| - 1).
\end{aligned}$$

This shows that

$$\begin{aligned}
X^v &\in \sigma(X^s)X^u + \pi X^u L(b, -qb|s| - 1) \\
&\subseteq \sigma(X^s)X^u + \pi L(b, -b|v| - 1),
\end{aligned}$$

where we used the relation $|v| = q|s| + |u|$. Thus, we can write

$$X^v = \sigma(X^s)X^u + \pi g_v, \quad g_v \in L(b, -b|v| - 1).$$

Now, we consider a general element

$$f = \sum_{v \in \mathbf{Z}_{\geq 0}^n} a_v X^v \in L(b, c), \quad \text{ord}_\pi a_v \geq b|v| + c.$$

By above, we can write

$$\begin{aligned}
f &= \sum_{0 \leq u_i < q} \sigma(\sum_v a_v X^{[v/q]})X^u + \pi \sum_v a_v g_v \\
&= \sum_{0 \leq u_i < q} \sigma(\sum_v a_v X^{[v/q]})X^u + \pi g,
\end{aligned}$$

where $[v/q] = ([v_1/q], \cdots, [v_n/q])$ denotes the integer part of the vector $v/q$ and

$$a_v g_v \in L(b, -b|v| - 1 + b|v| + c) = L(b, c - 1).$$

Thus, the limit

$$\lim_{|v| \to \infty} a_v g_v = 0$$

holds in $A_0$ and the sum

$$\pi g = \sum_v \pi a_v g_v \in L(b, c)$$

converges in $L(b, c)$ because this subspace of $A_0$ is closed. Since

$$b|v| + c \geq qb|[v/q]| + c,$$

we deduce

$$\sum_v a_v X^{[v/q]} \in L(qb, c).$$

Thus, we can write

$$f = \sum_{0 \leq u_i < q} \sigma(h_u)X^u + \pi g,$$

where
$$\pi g \in L(b,c), \quad h_u \in L(qb,c).$$
Iterating the above construction, we can then write
$$f = \sum_{0 \le u_i < q} \sigma(\sum_{j=0}^{\infty} \pi^j h_{uj}) X^u,$$
where
$$\pi^j h_{uj} \in L(qb,c), \quad \sum_{j=0}^{\infty} \pi^j h_{uj} \in L(qb,c).$$
This proves that
$$\Theta_u(f) = \sum_{j=0}^{\infty} \pi^j h_{uj} \in L(qb,c).$$
The proof is complete.                                                      □

For an arbitrary overconvergent Dwork operator on $A_0$, the situation can be reduced to the case of basic Dwork operators. We have

**Theorem 4.6.** *Let $\Theta$ be any fixed overconvergent Dwork operator on $A_0$. Choose rational numbers $b$ and $c_1$ with $0 < b \le b_\sigma$ such that for all $0 \le u_i < q$,*
$$\Theta(X^u) \in L(qb,c_1).$$
*Then, for all real number $c$, we have*
$$\Theta(L(b,c)) \subset L(qb, c+c_1).$$

*Proof.* For $f \in L(b,c)$, write
$$f = \sum_{0 \le u_i < q} \Theta_u(f)^\sigma X^u.$$
By the $\sigma^{-1}$-linearity of a Dwork operator, we deduce
$$\Theta(f) = \sum_{0 \le u_i < q} \Theta_u(f) \cdot \Theta(X^u).$$
Since
$$\Theta_u(f) \in L(qb,c), \quad \Theta(X^u) \in L(qb,c_1),$$
we conclude that
$$\Theta(f) \in L(qb, c+c_1).$$
The proof is complete.                                                      □

Thus, an overconvergent Dwork operator $\Theta$ on $A_0$ improves the rate of overconvergence. It takes an overconvergent element $f \in L(b,c)$ to a more overconvergent element $\Theta(f) \in L(qb, c+c_1)$. Theorem 4.6 shows that
$$\Theta(L(b,0) \otimes K) \subset L(qb,0) \otimes K \hookrightarrow L(b,0) \otimes K.$$
The last inclusion implies that $\Theta$ induces a completely continuous operator on the Banach space $L(b,0) \otimes K$, that is, the matrix of $\Theta$ is completely continuous in the sense of Serre [11]. Thus, the Fredholm determinant $\det(I - T\Theta)$ is well defined and $p$-adic entire. It is independent of the choice of $b$ satisfying $0 < b \le b_\sigma$. On the other hand, the Frobenius map $\sigma$ on $A_0$ is far from being completely continuous. It decreases the rate of overconvergence. One has

$$\sigma : L(b,c) \longrightarrow L(b/q,c).$$

For our applications, we need to generalize the above construction from the trivial rank one overconvergent $\sigma$-module $(A_0, \sigma)$ to an infinite rank nuclear overconvergent $\sigma$-module $(M, \phi)$.

Let $M$ be the Banach $A_0$-module with the formal basis $\vec{e} = \{e_1, e_2, \cdots\}$. That is, we have

$$M = \{\sum_j a_j e_j | a_j \in A_0\}.$$

Let

$$M^* = \mathrm{Hom}_{A_0}^{\mathrm{cont}}(M, \Omega^n A_0)$$

be the set of continuous $A_0$-linear homomorphisms from $M$ to $\Omega^n A_0$, where $\Omega^n A_0$ is the module of differential $n$-forms on $A_0$. As an abstract $A_0$-module, $\Omega^n A_0$ is isomorphic to $A_0$. Thus, $M^*$ is isomorphic to the continuous dual $M^\vee$ defined in section 2. However, we shall us $M^*$ later in the Monsky trace formula. The module $M^*$ is the Banach $A_0$-module with the orthonormal dual basis $\vec{e}^* = \{e_1^*, e_2^*, \cdots\}$, where

$$e_{j_1}^*(e_{j_2}) = \begin{cases} dX_1 \wedge \cdots \wedge dX_n, & \text{if } j_1 = j_2, \\ 0, & \text{if } j_1 \neq j_2. \end{cases}$$

That is, we have

$$M^* = \{\sum_j a_j e_j^* | a_j \in A_0, \ \lim_j \|a_j\| = 0\}.$$

The spaces $M$ and $M^*$ are too big for us. We often need to work with certain smaller subspaces when handling overconvergent $\sigma$-modules. For this purpose, we define the overconvergent subspaces of $M$ and $M^*$ by

$$M^\dagger = \{\sum_j a_j e_j | a_j \in A\},$$

$$M^{*\dagger} = \{\sum_j a_j e_j^* | a_j \in A, \ \lim_j \|a_j\| = 0\}.$$

These overconvergent spaces are, however, not complete. They are still a little big. What we also need to work on are some even smaller but complete subspaces of $M^\dagger$ and $M^{*\dagger}$, with explicit growth conditions. These spaces are defined by

$$M(b,c) = \{\sum_j a_j e_j \in M^\dagger | a_j \in L(b,c)\},$$

$$M^*(b,c) = \{\sum_j a_j e_j^* \in M^{*\dagger} | a_j \in L(b,c), \ \lim_j \|a_j\| = 0\},$$

where $b > 0$ and $c$ are rational numbers. Like $L(b,c)$, the spaces $M(b,c)$ and $M^*(b,c)$ are also complete with respect to the $\pi$-adic topology. One checks that $M(b,c)$ and $M^*(b,c)$ are $L(b,0)$-modules. More generally, we have the relations

$$L(b,c_1)M(b,c) \subset M(b,c_1+c),$$

$$L(b,c_1)M^*(b,c) \subset M^*(b,c_1+c).$$

It is easy to check that $M^*(b,0)$ is isomorphic to the continuous $L(b,0)$-linear dual of $M(b,0)$:

$$M^*(b,0) = \text{Hom}^{\text{cont}}_{L(b,0)}(M(b,0), L(b,0)).$$

Conversely, $M(b,0)$ is isomorphic to the continuous $L(b,0)$-linear dual of $M^*(b,0)$:

$$M(b,0) = \text{Hom}^{\text{cont}}_{L(b,0)}(M^*(b,0), L(b,0)).$$

The subspace $M(b,0)$ can be used to give an equivalent reformulation of the overconvergence of $\phi$ for a nuclear $\sigma$-module $(M, \phi)$. In fact, the following result is immediate from our definitions.

**Proposition 4.7.** *Let $(M, \phi)$ be a nuclear $\sigma$-module. Then, $\phi$ is overconvergent with respect to the formal basis $\vec{e}$ if and only if there are rational numbers $b > 0$ and $c$ such that all the entries of the matrix $B(X)$ of $\phi$ are in the space $L(b,c)$. Equivalently, we have $\phi(e_j) \in M(b,c)$ for all $j$.*

Note that the map $\phi$ itself is still not stable on $M(b,c)$ or $M(b,c) \otimes K$. The problem is caused by the $\sigma$-linearity acting on the coefficients. We have only the weaker inclusion

$$\phi(M(b,c) \otimes K) \subset M(b/q, c) \otimes K.$$

This is because

$$\sigma(L(b,c) \otimes K) \subset L(b/q, c) \otimes K.$$

A similar result holds for a family.

**Proposition 4.8.** *Let $(M, \phi_k)$ be a family of nuclear $\sigma$-modules parametrized by $k$ with the same underlying module $M$. Then, the family $\phi_k$ is uniformly overconvergent with respect to the formal basis $\vec{e}$ if and only if there are rational numbers $b > 0$ and $c$ independent of $k$ such that $\phi_k(e_j) \subset M(b,c)$ for all $j$ and all $k$.*

We now return to the definition of Dwork operators in the infinite rank setting. The extended Monsky trace formula expresses the L-function $L(\phi, T)$ of a nuclear overconvergent $\sigma$-module $\phi$ in terms of certain infinite rank Dwork operators acting on $M^*$. These Dwork operators are thus important in our investigation.

**Definition 4.9.** A Dwork operator $\Theta$ on $M^*$ is a continuous $R$-linear endomorphism of $M^*$ which is $\sigma^{-1}$-linear in the sense that

$$\Theta(\sigma(a)f) = a\Theta(f), \ a \in A_0, \ f \in M^*.$$

The Dwork operator $\Theta$ is called overconvergent if $\Theta$ is stable on the overconvergent subspace $M^{*\dagger}$, that is,

$$\Theta(M^{*\dagger}) \subset M^{*\dagger}.$$

The Dwork operator $\Theta$ is called nuclear if it satisfies the following nuclear condition: If we write

$$\Theta(X^v e^*_{j_2}) = \sum_{u, j_1} G^*_{\{u, j_1\}, \{v, j_2\}} X^u e^*_{j_1}, \ \ G^*_{\{u, j_1\}, \{v, j_2\}} \in R,$$

then for each fixed $u$ and each fixed $v$, we have

$$\lim_{j_1 \to \infty} \inf_{j_2} \text{ord}_\pi G^*_{\{u, j_1\}, \{v, j_2\}} = \infty.$$

If $M^*$ is of finite rank, then there are only finitely many $j_1$ and thus the nuclear condition in Definition 4.9 is automatically satisfied. In general, the nuclear condition can be viewed as a uniform condition. The $\sigma^{-1}$-linearity and the continuity of $\Theta$ imply the infinite $\sigma^{-1}$-linearity of $\Theta$, namely,

$$\Theta(\sum_j \sigma(a_j)\alpha_j e_j^*) = \sum_j a_j\Theta(\alpha_j e_j^*)$$

for all elements $a_j \in A_0$ with $\lim_j \|a_j\| = 0$ and all elements $\alpha_j \in A_0$. The norm of $\Theta$ is defined to be

$$\|\Theta\| = p^{-\mathrm{ord}_\pi\Theta},$$

where $\mathrm{ord}_\pi\Theta$ is the smallest non-negative integer $i$ such that

$$\Theta(M^*) \subset \pi^i M^*.$$

Since $\Theta$ is integral, we have $\|\Theta\| \leq 1$.

**Definition 4.10.** A Dwork operator $\Theta$ on $M^*$ is called **contracting** if there are rational numbers $b > 0$, $c$ and $c_1$ such that

$$\Theta(M^*(b,c)) \subset M^*(qb, c+c_1).$$

Equivalently, for all $j$, we have

$$\Theta(L(b,c)e_j^*) \subset M^*(qb, c+c_1).$$

It is easy to check that a contracting Dwork operator is automatically overconvergent. Conversely, if $M^*$ is of finite rank, then any overconvergent Dwork operator on $M^*$ is automatically contracting by a simple application of Theorem 4.6. In the infinite rank case, the contracting condition in Definition 4.10 can be viewed as another uniform condition. It is stronger than the overconvergent condition in Definition 4.9. We now show that the Fredholm determinant $\det(I - T\Theta)$ is well defined and $p$-adic entire for a nuclear contracting Dwork operator $\Theta$ on $M^*$.

Let $\Theta$ be a nuclear contracting Dwork operator on $M^*$ with $0 < b \leq b_\sigma$, $c$ and $c_1$ satisfying the requirements of Definition 4.10. For simplicity of notations, we shall now assume that $b$ is a positive integer. As indicated before, this can always be achieved by going to a totally ramified finite extension of $R$ if necessary. The space $M^*(b,c) \otimes K$ is a $K$-Banach space with the topological basis $\{\pi^{b\lfloor u\rfloor}X^u e_j^*\}$. The contracting condition

$$\Theta(M^*(b,c)) \subset M^*(qb, c+c_1)$$

shows that the map $\Theta$ induces a continuous $K$-linear endomorphism of the Banach space $M^*(b,c) \otimes K$. We now show that the matrix of $\Theta$ is completely continuous in the sense of Serre.

Let $G = (G_{\{u,j_1\},\{v,j_2\}})$ be the infinite matrix of the $K$-linear map $\Theta$ acting on $M^*(b,c) \otimes K$ defined with respect to the topological basis $\{\pi^{b\lfloor u\rfloor}X^u e_j^*\}$. The subscripts $\{u,j_1\}$ denote the row index of $G$. The subscripts $\{v,j_2\}$ denote the column index of $G$. By our definition of the matrix $G$, we have

$$\Theta(\pi^{b\lfloor v\rfloor}X^v e_{j_2}^*) = \sum_{u,j_1} G_{\{u,j_1\},\{v,j_2\}}\pi^{b\lfloor u\rfloor}X^u e_{j_1}^*.$$

In terms of the notation in Definition 4.9, we have

$$G_{\{u,j_1\},\{v,j_2\}} = \pi^{b\lfloor v\rfloor - b\lfloor u\rfloor}G_{\{u,j_1\},\{v,j_2\}}^*.$$

We want to show that the matrix $G$ is completely continuous. Namely,

**Lemma 4.11.** *Let $\Theta$ be a nuclear contracting Dwork operator on $M^*$. For any constant $C > 0$, except for finitely many row indices $\{u, j_1\}$, we always have the bound*

$$\mathrm{ord}_\pi G_{\{u, j_1\}, \{v, j_2\}} \geq C.$$

*Proof.* From the contracting condition (one can take $c = 0$)

$$\Theta(M^*(b, 0)) \subset M^*(qb, c_1),$$

we deduce that

$$
\begin{aligned}
\mathrm{ord}_\pi G_{\{u, j_1\}, \{v, j_2\}} &\geq qb|u| + c_1 - b|u| \\
&\geq (q-1)b|u| + c_1.
\end{aligned}
$$

The right side is uniformly large as long as $|u|$ is large. Thus, we may restrict our attention to the finitely many $u$ with small $|u|$. Since $\Theta$ is integral (mapping $M^*$ to $M^*$), we deduce that

$$
\begin{aligned}
\mathrm{ord}_\pi G_{\{u, j_1\}, \{v, j_2\}} &= b(|v| - |u|) + \mathrm{ord}_\pi G^*_{\{u, j_1\}, \{v, j_2\}} \\
&\geq b(|v| - |u|).
\end{aligned}
$$

The right side is again uniformly large for fixed $|u|$ as long as $|v|$ is large. Thus, we may assume that both $|u|$ and $|v|$ are small or fixed. Now, the nuclear condition of $\Theta$ is

$$\lim_{j_1 \to \infty} \inf_{j_2} \mathrm{ord}_\pi G^*_{\{u, j_1\}, \{v, j_2\}} = \infty.$$

It implies that

$$\lim_{j_1 \to \infty} \inf_{j_2} \mathrm{ord}_\pi G_{\{u, j_1\}, \{v, j_2\}} = \infty.$$

Thus,

$$\mathrm{ord}_\pi G_{\{u, j_1\}, \{v, j_2\}} \geq C$$

uniformly for fixed $|u|$, fixed $|v|$, all $j_2$ and all large $j_1$. There are thus only finitely many row indices $\{u, j_1\}$ which may fail to satisfy the above inequality. The proof is complete. $\qquad\square$

Lemma 4.11 shows that the the Fredholm determinant

$$\det(I - T\Theta | M^*(b, 0) \otimes K) = \det(I - TG)$$

is a $p$-adic entire function. This determinant is independent of our choice of sufficiently small positive number $b$ as the matrices $G$ for different $b$ are conjugate.

**Definition 4.12.** Let $\Theta$ be a nuclear contracting Dwork operator on $M^*$. Define the Fredholm determinant $\det(I - T\Theta | M^*)$ of $\Theta$ on $M^*$ to be the following entire function

$$\det(I - T\Theta | M^*) = \det(I - T\Theta | M^*(b, 0) \otimes K) = \det(I - TG),$$

where $b$ is any positive number satisfying $0 < b \leq b_\sigma$ and the requirements of Definition 4.10.

Instead of using the matrix $G = (G_{\{u, j_1\}, \{v, j_2\}})$ in the definition of the Fredholm determinant, we could also use the matrix $G^* = (G^*_{\{u, j_1\}, \{v, j_2\}})$ from Definition 4.9. The result is the same. That is,

$$\det(I - TG^*) = \det(I - TG)$$

The advantage for using $G$ is to be able to see the completely continuous shape of the matrix $G$ right away. It makes later explicit estimates of Newton polygon more transparent. It has, however, one drawback. Namely, the entries in $G$ are not integral although they are indeed bounded. On the other hand, the entries in $G^*$ are always integral. This is useful to prove explicit continuity result. We state a simple result here, whose proof follows immediately from the use of the integral matrix $G^*$ in the Fredholm determinant.

**Lemma 4.13.** *Let $\Theta_1$ and $\Theta_2$ be two nuclear contracting Dwork operators on $M^*$. Assume that*

$$\Theta_1 \equiv \Theta_2 \ (\mathrm{mod}\,\pi^i)$$

*for some non-negative integer $i$. Then, we have the congruence*

$$\det(I - T\Theta_1|M^*) \equiv \det(I - T\Theta_2|M^*) \ (\mathrm{mod}\,\pi^i).$$

We now turn to the family version of the above results.

**Definition 4.14.** Let $\Theta_k$ be a family of nuclear contracting Dwork operators on the same underlying space $M^*$, parameterized by $k$. We say that the family $\Theta_k$ is **uniform** if the following two conditions hold. First, the $\Theta_k$ satisfy the uniform nuclear condition: if we write

$$\Theta_k(X^v e_{j_2}^*) = \sum_{u,j_1} G_{\{u,j_1\},\{v,j_2\}}^*(k) X^u e_{j_1}^*,$$

then for fixed $u$ and fixed $v$, we have

$$\lim_{j_1 \to \infty} \inf_{j_2,k} \mathrm{ord}_\pi G_{\{u,j_1\},\{v,j_2\}}^*(k) = \infty.$$

Second, the $\Theta_k$ satisfy the uniform contracting condition

$$\Theta_k(M^*(b,c)) \subset M^*(qb, c+c_1),$$

where the constants $b > 0$, $c$ and $c_1$ are independent of the parameter $k$.

**Definition 4.15.** We say that the family $\Theta_k$ is **uniformly continuous** in $k$ if $\|\Theta_{k_1} - \Theta_{k_2}\|$ is uniformly small whenever $k_1$ is close to $k_2$.

Let

$$G(k) = (G_{\{u,j_1\},\{v,j_2\}}(k))$$

be the matrix of $\Theta_k$ with respect to the topological basis $\{\pi^{b|v|} X^u e_j^*\}$ of the Banach space $M^*(b,c) \otimes K$. Then, we have

**Lemma 4.16.** *Let $\Theta_k$ be a uniform family of nuclear contracting Dwork operators on $M^*$. For any constant $C > 0$, except for finitely many fixed row indices $\{u,j_1\}$ depending on $C$, we have the bound*

$$\mathrm{ord}_\pi G_{\{u,j_1\},\{v,j_2\}}(k) \geq C$$

*uniformly for all $k$.*

The proof is the same as the proof of Lemma 4.11. Our uniform condition on $\Theta_k$ guarantees that the proof goes through in a uniform way.

**Corollary 4.17.** *Let $\Theta_k$ be a uniform family of nuclear contracting Dwork operators on $M^*$. Then the family of Fredholm determinants*

$$\det(I - T\Theta_k|M^*) = \det(I - TG(k))$$

*is a family of uniformly entire functions. If, in addition, the family $\Theta_k$ is uniformly continuous in $k$, then $\det(I - T\Theta_k|M^*)$ is a uniformly continuous family of uniformly entire functions and we have*

$$\| \det(I - T\Theta_{k_1}|M^*) - \det(I - T\Theta_{k_2}|M^*) \| \leq \|\Theta_{k_1} - \Theta_{k_2}\|.$$

## 5. THE EXTENDED MONSKY TRACE FORMULA

In this section, we extend the Monsky trace formula from finite rank case to infinite rank case. One way to achieve this would be to extend the whole proof and setup of the original Monsky trace formula to infinite rank setting. This should be possible in principal but would be very long. Instead of reproving it in a more complicated situation, we shall take a much shorter limiting approach which simply uses the finite rank version of the Monsky trace formula. The reader is referred to [17] for an explicit description and a proof of the Monsky trace formula in the finite rank case. Here we give a summary of the description in the affine $n$-space case $\mathbf{A}^n$.

For each integer $0 \leq i \leq n$, let $\Omega^i A$ be the finite free $A$-module of differential $i$-forms on $A$ over $R$. In particular, $\Omega^n A$ is the free rank one $A$-module with the generator $dX_1 \wedge \cdots \wedge dX_n$. The Frobenius map $\sigma$ extends to an injective $R$-linear endomorphism $\sigma_i$ of $\Omega^i A$. It satisfies

$$\sigma_i(af) = \sigma(a)\sigma_i(f), \ a \in A, \ f \in \Omega^i A.$$

Thus, the pair $(\Omega^i A, \sigma_i)$ becomes a finite rank overconvergent $\sigma$-module. There is a trace map

$$\mathrm{Tr}_i : \Omega^i A \longrightarrow \sigma_i(\Omega^i A)$$

such that for $a \in A, f \in \Omega^i A$, we have

$$\mathrm{Tr}_i(\sigma(a)f) = \sigma(a)\mathrm{Tr}_i(f).$$

Furthermore, as endomorphisms of $\Omega^i A$, we have

(5.1)                    $$\sigma_i^{-1} \circ \mathrm{Tr}_i \circ \sigma_i = [\bar{A} : \sigma(\bar{A})] = q^n,$$

where $\bar{A}$ means the reduction modulo $\pi$ of $A$. Thus, the map $\sigma_i$ acting on $\Omega^i A$ has a one-sided left inverse when tensored with $\mathbf{Q}$. The map $\sigma_i^{-1} \circ \mathrm{Tr}_i$ acting on $\Omega^i A$ is $\sigma^{-1}$-linear in the sense that

$$(\sigma_i^{-1} \circ \mathrm{Tr}_i)(\sigma(a)f) = a(\sigma_i^{-1} \circ \mathrm{Tr}_i)(f), \ a \in A, \ f \in \Omega^i A.$$

Both $\sigma_i$ and the trace map $\mathrm{Tr}_i$ extend to $\Omega^i A_0 = \Omega A \otimes A_0$ by continuity. In particular, the map $\sigma_i^{-1} \circ \mathrm{Tr}_i$ acting on $\Omega^i A_0$ is an overconvergent Dwork operator. This construction corresponds to the trivial rank one overconvergent $\sigma$-module $(A_0, \sigma)$. We need to extend it to an infinite rank overconvergent nuclear $\sigma$-module.

Let $(M, \phi)$ be an overconvergent nuclear $\sigma$-module with the formal basis $\vec{e}$. Recall that in section 4, we defined a new dual $M^*$ by

$$M^* = \mathrm{Hom}_{A_0}^{\mathrm{cont}}(M, \Omega^n A_0) = \{\sum_j a_j e_j^* | a_j \in A_0, \ \lim_j \|a_j\| = 0\}.$$

where $e^* = \{e_1^*, e_2^*, \cdots, \}$ denotes the dual basis of $\vec{e}$ defined by:

$$e_{j_1}^*(e_{j_2}) = \begin{cases} dX_1 \wedge \cdots \wedge dX_n, & \text{if } j_1 = j_2, \\ 0, & \text{if } j_1 \neq j_2. \end{cases}$$

The $A_0$-linear map

$$\phi : M^\sigma \longrightarrow M$$

induces an $A_0$-linear map on their duals:

$$\phi^* : M^* \longrightarrow M^{*\sigma} = M^{\sigma*}, \ \ \phi^*(f) = f \circ \phi, \ f \in M^*,$$

where

$$M^{\sigma*} = \{\sum a_j \otimes e_j^* | a_j \in A_0, \ \lim_j \|a_j\| = 0\}, \ \sigma(a_j) \otimes e_j^* = a_j e_j^*.\}$$

The matrix of $\phi^*$ is just the transpose of the matrix of $\phi$.

**Definition 5.1.** Let $(M, \phi)$ be an overconvergent nuclear $\sigma$-module. The associated Dwork operator $\Theta$ on $M^*$ is the $R$-linear endomorphism defined by

$$\Theta(f)(m) = (\sigma_n^{-1} \circ \mathrm{Tr}_n)(f(\phi(m))), \ f \in M^*, \ m \in M,$$

where $\phi$ is viewed as a $\sigma$-linear map from $M$ to itself.

One checks that $\Theta(f)$ is in $M^*$. Furthermore, the map $\Theta$ is $\sigma^{-1}$-linear. That is, for $a \in A_0$ and $f \in M^*$, we have

$$\Theta(\sigma(a)f) = a\Theta(f).$$

This map $\Theta$ is nuclear and contracting as we shall show now.

**Lemma 5.2.** *Let $(M, \phi)$ be an overconvergent nuclear $\sigma$-module. Then, the associated Dwork operator $\Theta$ on $M^*$ is a nuclear contracting Dwork operator. In particular, the Fredholm determinant $\det(I - T\Theta | M^*)$ is well defined and $p$-adically entire.*

*Proof.* Let $B = (B_{j_1, j_2})$ be the matrix of $\phi$ with respect to the formal basis $\vec{e}$, where the entries $B_{j_1, j_2}$ are in $A$. That is,

$$\phi(e_{j_2}) = \sum_{j_1} B_{j_1, j_2} e_{j_1}.$$

Let

$$m = \sum_j a_j e_j, \ a_j \in A_0$$

be a typical element of $M$. For $a \in A_0$, one checks that

$$
\begin{aligned}
\Theta(ae_j^*)(m) &= (\sigma_n^{-1} \circ \mathrm{Tr}_n)(ae_j^*(\phi(m))) \\
&= (\sigma_n^{-1} \circ \mathrm{Tr}_n)(ae_j^*(\sum_{j_2} \sigma(a_{j_2})\phi(e_{j_2}))) \\
&= (\sigma_n^{-1} \circ \mathrm{Tr}_n)(ae_j^*(\sum_{j_2, j_1} \sigma(a_{j_2}) B_{j_1, j_2} e_{j_1})) \\
&= \sum_{j_2} a_{j_2}(\sigma_n^{-1} \circ \mathrm{Tr}_n)(aB_{j, j_2} dX_1 \wedge \cdots \wedge dX_n) \\
&= \sum_{j_2} D_{j, j_2}(a) a_{j_2} dX_1 \wedge \cdots \wedge dX_n \\
&= \sum_{j_2} D_{j, j_2}(a) e_{j_2}^*(m),
\end{aligned}
$$

where the elements $D_{j,j_2}(a) \in A_0$ are defined by

$$D_{j,j_2}(a) = \frac{(\sigma_n^{-1} \circ \mathrm{Tr}_n)(a B_{j,j_2} dX_1 \wedge \cdots \wedge dX_n)}{dX_1 \wedge \cdots \wedge dX_n}.$$

Changing the notation, we obtain the relation in $M^*$:

$$\Theta(a e_{j_1}^*) = \sum_{j_2} D_{j_1,j_2}(a) e_{j_2}^*.$$

To prove that $\Theta$ is indeed a nuclear Dwork operator on $M^*$, it suffices to check

$$\lim_{j_2 \to \infty} \inf_{j_1} \mathrm{ord}_\pi D_{j_1,j_2}(X^v) = \infty.$$

But this holds since

$$\lim_{j_2 \to \infty} \inf_{j_1} \mathrm{ord}_\pi(X^v B_{j_1,j_2}) = \lim_{j_2 \to \infty} \inf_{j_1} \mathrm{ord}_\pi(B_{j_1,j_2}) = \infty.$$

This last limit follows from the condition that $\phi$ is a nuclear $\sigma$-module. Thus, the nuclear condition of $\Theta$ is satisfied.

We now prove that $\Theta$ is contracting. For this purpose, we choose rational numbers $b > 0$ and $c$ such that

$$B_{j_1,j_2} \in L(b,c)$$

uniformly for all $j_1$ and $j_2$. This is possible since $\phi$ is overconvergent. Choosing smaller $b$ and larger $c$ if necessary, we may assume that

$$(\sigma_n^{-1} \circ \mathrm{Tr}_n)(L(b,2c) dX_1 \wedge \cdots \wedge dX_n) \subset L(qb, c + c_1) dX_1 \wedge \cdots \wedge dX_n.$$

This is again possible since $\sigma_n^{-1} \circ \mathrm{Tr}_n$ is an overconvergent Dwork operator on the rank one $A_0$-module $\Omega^n A_0$ and we can use Theorem 4.6. With such a choice of $b > 0$, $c$ and $c_1$, we deduce that

$$D_{j_2,j_1}(L(b,c)) \subset L(qb, c + c_1)$$

uniformly for all $j_1$ and $j_2$. This gives the desired contracting property

$$\Theta(M^*(b,c)) \subset M^*(qb, c + c_1).$$

The proof of Lemma 5.2 is complete. $\qquad\blacksquare$

For later applications, we derive an explicit estimate for the entries of the matrix of the Dwork operator $\Theta$ in a special case.

**Lemma 5.3.** *Assume that there is a sequence of non-negative integers $o(j)$ such that*

$$\phi(e_j) \in \pi^{o(j)} M(b, c_1)$$

*for all $j$. Then, there is a constant $c$ independent of $j_1$ and $j_2$ such that*

$$D_{j_1,j_2}(L(b,0)) \subset \pi^{o(j_2)} L(qb, c).$$

*Proof.* By our assumption, we can write

$$B_{j_1,j_2} = \pi^{o(j_2)} C_{j_1,j_2}, \ \ C_{j_1,j_2} \in L(b, c_1).$$

Then, for $a \in L(b, 0)$, we deduce from Theorem 4.6 the following estimate:

$$D_{j_1,j_2}(a) = \pi^{o(j_2)} \frac{(\sigma_n^{-1} \circ \mathrm{Tr}_n)(a C_{j_1,j_2} dX_1 \wedge \cdots \wedge dX_n)}{dX_1 \wedge \cdots \wedge dX_n} \in \pi^{o(j_2)} L(qb, c).$$

The lemma is proved. $\qquad\blacksquare$

**Lemma 5.4.** *Assume that there is a sequence of non-negative integers $o(j)$ such that*

$$\phi(e_j) \in \pi^{o(j)} M(b, c_1)$$

*for all $j$. Write*

$$\Theta(\pi^{b|v|} X^v e_{j_2}^*) = \sum_{u, j_1} G_{\{u, j_1\}, \{v, j_2\}} \pi^{b|u|} X^u e_{j_1}^*.$$

*Then, there is a constant $c$ such that*

$$\mathrm{ord}_\pi G_{\{u, j_1\}, \{v, j_2\}} \geq (q-1)b|u| + o(j_1) + c.$$

*Proof.* Taking $a = \pi^{b|v|} X^v \in L(b, 0)$, we obtain from Lemma 5.3 that

$$D_{j_2, j_1}(\pi^{b|v|} X^v) \in \pi^{o(j_1)} L(qb, c).$$

But we also have the relation

$$\Theta(\pi^{b|v|} X^v e_{j_2}^*) = \sum_{j_1} D_{j_2, j_1}(\pi^{b|v|} X^v) e_{j_1}^*.$$

Thus,

$$\sum_u G_{\{u, j_1\}, \{v, j_2\}} \pi^{b|u|} X^u = D_{j_2, j_1}(\pi^{b|v|} X^v) \in \pi^{o(j_1)} L(qb, c).$$

It follows that

$$\begin{aligned}
\mathrm{ord}_\pi G_{\{u, j_1\}, \{v, j_2\}} &\geq o(j_1) + qb|u| + c - b|u| \\
&\geq (q-1)b|u| + o(j_1) + c.
\end{aligned}$$

The Lemma is proved. $\qquad\blacksquare$

In our later applications, the sequence $o(j)$ will go to infinity when $j$ goes to infinity.

To get the extended Monsky trace formula, we have to generalize the above results to differential $i$-forms. For each non-negative integer $0 \leq i \leq n$, let

$$\Omega^i M = M \otimes_{A_0} \Omega^i A_0 = \{\sum_j a_j e_j | a_j \in \Omega^i A_0\},$$

where the tensor product denotes the formal tensor product. This is again a Banach $A_0$-module with the formal basis

$$\vec{e} \otimes \wedge^i = \{e_j \otimes (dX_{\ell_1} \wedge \cdots \wedge dX_{\ell_i}) | j \geq 1, 1 \leq \ell_1 < \ell_2 < \cdots < \ell_i \leq n\}.$$

That is,

$$\Omega^i M = \{\sum_j \sum_{1 \leq \ell_1 < \ell_2 < \cdots < \ell_i \leq n} a_j(\ell_1, \cdots, \ell_i) e_j(\ell_1, \cdots, \ell_i) \mid a_j(\ell_1, \cdots, \ell_i) \in A_0\},$$

where

$$e_j(\ell_1, \cdots, \ell_i) = e_j \otimes (dX_{\ell_1} \wedge \cdots \wedge dX_{\ell_i}).$$

Let

$$\phi_i = \phi \otimes \sigma_i,$$

where $(\Omega^i A_0, \sigma_i)$ is the previously mentioned finite rank overconvergent $\sigma$-module given by the action of $\sigma_i$ on $\Omega^i A_0$. Then, the pair $(\Omega^i M, \phi_i)$ becomes an overconvergent nuclear $\sigma$-module. As above, we define a new dual $A_0$-module by

$$
\begin{aligned}
M_i^* &= \operatorname{Hom}_{A_0}^{\mathrm{cont}}(\Omega^i M, \Omega^n A_0) \\
&= \{\sum_j a_j e_j^* | a_j \in \operatorname{Hom}_{A_0}(\Omega^i A_0, \Omega^n A_0), \ \lim_j \|a_j\| = 0\}.
\end{aligned}
$$

This is the set of continuous $A_0$-linear maps from $\Omega^i M$ to $\Omega^n A_0$. It is a Banach $A_0$-module with the orthonormal basis

$$
\bar{e}^* \otimes \wedge^i = \{e_j^* \otimes (dX_{\ell_1} \wedge \cdots \wedge dX_{\ell_i}) | j \geq 1, 1 \leq \ell_1 < \ell_2 < \cdots < \ell_i \leq n\}.
$$

That is,

$$
M_i^* = \{\sum_j \sum_{1 \leq \ell_1 < \ell_2 < \cdots < \ell_i \leq n} a_j(\ell_1, \cdots, \ell_i) e_j^*(\ell_1, \cdots, \ell_i) | a_j \in A_0, \ \lim_j \|a_j\| = 0\},
$$

where

$$
e_j^*(\ell_1, \cdots, \ell_i) = e_j^* \otimes (dX_{\ell_1} \wedge \cdots \wedge dX_{\ell_i}).
$$

The overconvergent subspace of $M_i^*$ is defined to be

$$
M_i^{*\dagger} = \{\sum_j a_j e_j^* | a_j \in \operatorname{Hom}_A(\Omega^i A, \Omega^n A), \ \lim_j \|a_j\| = 0\}.
$$

In terms of the basis $\bar{e}^* \wedge^i$, we have the description

$$
M_i^{*\dagger} = \{\sum_j \sum_{1 \leq \ell_1 < \ell_2 < \cdots < \ell_i \leq n} a_j(\ell_1, \cdots, \ell_i) e_j^*(\ell_1, \cdots, \ell_i) | a_j \in A, \ \lim_j \|a_j\| = 0\}.
$$

For rational numbers $b > 0$ and $c$, we similarly define

$$
M_i^*(b, c) = \{\sum_j \sum_{\ell_1, \cdots, \ell_i} a_j(\ell_1, \cdots, \ell_i) e_j^*(\ell_1, \cdots, \ell_i) | a_j \in L(b, c), \ \lim_j \|a_j\| = 0\},
$$

where the $\ell_j$ satisfy

$$
1 \leq \ell_1 < \ell_2 < \cdots < \ell_i \leq n.
$$

The module $M_i^*$ (resp. $M_i^*(b, c)$) is just $\binom{n}{i}$ copies of $M^*$ (resp. $M^*(b, c)$).

**Definition 5.5.** For each integer $0 \leq i \leq n$, the $i$-th Dwork operator $\Theta_i$ on $M_i^*$ associated to $\phi$ is defined by

$$
\Theta_i(f)(m) = (\sigma_n^{-1} \circ \operatorname{Tr}_n)(f(\phi_i(m))), \ m \in \Omega^i M, \ f \in M_i^*.
$$

This is again a $\sigma^{-1}$-linear operator.

Similarly to Lemma 5.2, we have

**Lemma 5.6.** *Let $(M, \phi)$ be an overconvergent nuclear $\sigma$-modules. Then for each integer $0 \leq i \leq n$, the $i$-th associated Dwork operator $\Theta_i$ on $M_i^*$ is nuclear and contracting. In particular, the Fredholm determinant $\det(I - T\Theta_i | M_i^*)$ is well defined and $p$-adically entire.*

Similar results hold for a family.

**Lemma 5.7.** *Let $(M, \phi(k))$ be a family of uniformly overconvergent nuclear $\sigma$-modules. Then for each $0 \leq i \leq n$, the $i$-th associated Dwork operator $\Theta_i(k)$ is a uniform family of nuclear contracting Dwork operators on $M_i^*$ parametrized by $k$. In particular, the Fredholm determinant $\det(I - T\Theta_i(k)|M_i^*)$ is a family of uniformly entire functions. If, in addition, the family $\phi(k)$ is uniformly continuous in $k$, then for each $0 \leq i \leq n$, the corresponding family of Dwork operators $\Theta_i(k)$ is also uniformly continuous in $k$:*

$$\|\Theta_i(k_1) - \Theta_i(k_2)\| \leq \|\phi(k_1) - \phi(k_2)\|.$$

*In this case, the Fredholm determinant $\det(I - T\Theta_i(k)|M_i^*)$ forms a family of uniformly continuous and uniformly entire functions. Furthermore,*

$$\|\det(I - T\Theta_i(k_1)|M_i^*) - \det(I - T\Theta_i(k_2)|M_i^*)\| \leq \|\phi(k_1) - \phi(k_2)\|.$$

The only thing that needs to be checked is the explicit continuity result. By Corollary 4.17, we only need to check that

$$\|\Theta_i(k_1) - \Theta_i(k_2)\| \leq \|\phi(k_1) - \phi(k_2)\|.$$

But this follows from our definition of the Dwork operator which is integral.

We can now state and prove the infinite rank version of the Monsky trace formula.

**Theorem 5.8.** *Let $(M, \phi)$ be a nuclear overconvergent $\sigma$-module. Then, we have the formula for the L-function:*

$$L(\phi, T) = \prod_{i=0}^{n} \det(I - T\Theta_i|M_i^*)^{(-1)^{i-1}}.$$

*In particular, the L-function $L(\phi, T)$ is a p-adic meromorphic function.*

*Proof.* If $\phi$ is of finite rank, this is proved in the appendix of [17]. We now show that it holds in infinite rank case as well, using a limiting argument.

For each positive integer $k$ including $\infty$, we define a new nuclear overconvergent $\sigma$-module $\phi_k$ on $M$ by

$$\phi_k(e_j) = \begin{cases} \phi(e_j), & \text{if } j \leq k, \\ 0, & \text{if } j > k. \end{cases}$$

Thus, $\phi = \phi_\infty$. Let the parameter $k$ vary in the set of positive integers including $\infty$ with the sequence topology. It is clear that the family $(M, \phi_k)$ is uniformly overconvergent and uniformly continuous. For $0 \leq i \leq n$, let $\Theta_i(k)$ be the $i$-th Dwork operator associated to $(M, \phi_k)$. Thus, Lemma 5.7 shows that for each $0 \leq i \leq n$, the Fredholm determinant $\det(I - T\Theta_i(k)|M_i^*)$ is a uniformly continuous family of uniformly entire functions. In particular,

$$\lim_{k \to \infty} \det(I - T\Theta_i(k)|M_i^*) = \det(I - T\Theta_i(\infty)|M_i^*)$$
$$= \det(I - T\Theta_i|M_i^*).$$

Let now $k$ be a finite positive integer. Let $M_k$ be the Banach $A_0$-submodule of $M$ defined by

$$M_k = \{\sum_{j > k} A_0 e_j\}.$$

It has the formal basis $\{e_{k+1}, e_{k+2}, \cdots, \}$. The map $\phi_k$ becomes the zero map on the submodule $M_k$ and hence stable on $M_k$. The quotient $(M/M_k, \phi_k)$ is a finite rank overconvergent $\sigma$-module, where

$$M/M_k \cong \{\sum_{j=1}^{k} A_0 e_j\}.$$

Denote the finite rank overconvergent $\sigma$-module $(M/M_k, \phi_k)$ by $(N_k, \psi_k)$. Thus, we have an exact sequence of nuclear overconvergent $\sigma$-modules:

$$0 \to (M_k, \phi_k) \to (M, \phi_k) \to (N_k, \psi_k) \to 0.$$

The family $(N_k, \psi_k)$ is uniformly overconvergent but not necessarily continuous in $k$ since the underlying module $N_k$ varies with $k$. The dual $A_0$-module of the finite rank $A_0$-module $\Omega^i N_k$ is given by

$$(N_k)_i^* = \mathrm{Hom}_{A_0}(\Omega^i N_k, \Omega^n A_0) = \{\sum_{j=1}^{k} \mathrm{Hom}_{A_0}(\Omega^i A_0, \Omega^n A_0) e_j^*\}.$$

Let $\Lambda_i(k)$ be the $i$-th Dwork operator associated to the finite rank overconvergent $\sigma$-module $(N_k, \psi_k)$. Since $\phi_k$ is the zero map on the submodule $M_k$ of $M$, we deduce that $\Theta_i(k)$ is the zero map on the submodule $M_{i,k}^*$ of $M_i^*$, where

$$M_{i,k}^* = \{\sum_{j>k} \mathrm{Hom}_{A_0}(\Omega^i A_0, \Omega^n A_0) e_j^*\}.$$

The restriction of $\Theta_i(k)$ to the quotient $M_i^*/M_{i,k}^* = (N_k)_i^*$ becomes the operator $\Lambda_i(k)$ on $(N_k)_i^*$. Namely, we have the exact sequence of Dwork operators:

$$0 \to ((N_k)_i^*, \Lambda_i(k)) \to (M_i^*, \Theta_i(k)) \to (M_{i,k}^*, 0) \to 0.$$

Thus implies that

$$\det(I - T\Theta_i(k)|M_i^*) = \det(I - T\Lambda_i(k)|(N_k)_i^*).$$

Now, since $\psi_k$ is of finite rank, we can apply the finite rank Monsky trace formula and deduce that for each finite $k$,

$$L(\psi_k, T) = \prod_{i=0}^{n} \det(I - T\Lambda_i(k)|(N_k)_i^*)^{(-1)^{i-1}}.$$

On the other hand, one checks directly from the Euler product definition of L-functions that we have

$$L(\phi_k, T) = L(\psi_k, T).$$

Since $\phi_k$ is uniformly continuous in $k$, the Euler product definition shows that

$$\lim_{k \to \infty} L(\phi_k, T) = L(\phi, T).$$

Combining with the previous equation, we deduce

$$\lim_{k \to \infty} L(\psi_k, T) = L(\phi, T).$$

Putting these together, we obtain

$$
\begin{aligned}
L(\phi, T) &= \lim_{k \to \infty} L(\psi_k, T) \\
&= \lim_{k \to \infty} \prod_{i=0}^{n} \det(I - T\Lambda_i(k)|(N_k)_i^*)^{(-1)^{i-1}} \\
&= \lim_{k \to \infty} \prod_{i=0}^{n} \det(I - T\Theta_i(k)|M_i^*)^{(-1)^{i-1}} \\
&= \prod_{i=0}^{n} \lim_{k \to \infty} \det(I - T\Theta_i(k)|M_i^*)^{(-1)^{i-1}} \\
&= \prod_{i=0}^{n} \det(I - T\Theta_i|M_i^*)^{(-1)^{i-1}}.
\end{aligned}
$$

The proof is complete. □

As an application, we obtain the following result for a family.

**Corollary 5.9.** *Let $\phi_k$ be a family of uniformly overconvergent nuclear $\sigma$-modules. Then the L-functions $L(\phi_k, T)$ form a family of uniformly meromorphic functions. If, in addition, the family $\phi_k$ is uniformly continuous in $k$, then the family of L-functions $L(\phi_k, T)$ parameterized by $k$ is a strong family of meromorphic functions.*

For later applications, we need to generalize the above result to the product of a suitable infinite families of L-functions. The precise result is as follows.

**Lemma 5.10.** *Assume that for each non-negative integer $i \geq 0$, we are given a family of uniformly continuous and uniformly overconvergent nuclear $\sigma$-modules $\phi_{k,i}$ parametrized by $k$. Let*

$$
L(k, T) = \prod_{i=0}^{\infty} L(\phi_{k,i}, \pi^i T).
$$

*Then, the family of functions $L(k, T)$ parameterized by $k$ is a strong family of meromorphic functions.*

*Proof.* By Lemma 5.7 and Corollary 5.9, we can write

$$
L(\phi_{k,i}, T) = \frac{f_{k,i}(T)}{g_{k,i}(T)},
$$

where for each fixed integer $i \geq 0$, both $f_{k,i}(T)$ and $g_{k,i}(T)$ form a family of uniformly continuous and uniformly entire functions parametrized by $k$. Furthermore, both $f_{k,i}(T)$ and $g_{k,i}(T)$ are power series with coefficients in $R$ and with constant term 1. They are a finite product of the Fredholm determinants as given in Lemma 5.7. Thus, by our definition of $L(k, T)$, we can write

$$
L(k, T) = \frac{f(k, T)}{g(k, T)},
$$

where

$$
f(k, T) = \prod_{i=0}^{\infty} f_{k,i}(\pi^i T)
$$

and

$$g(k,T) = \prod_{i=0}^{\infty} g_{k,i}(\pi^i T).$$

It suffices to prove that the family $f(k,T)$ (resp. the family $g(k,T)$) is a family of uniformly continuous and uniformly entire functions. We shall work with $f(k,T)$ as the proof for $g(k,T)$ is completely similar.

For a positive integer $j$, let

$$f_j(k,T) = \prod_{i=0}^{j-1} f_{k,i}(\pi^i T).$$

Since this is a finite product, we deduce that for each $j$, the family $f_j(k,T)$ parametrized by $k$ is a uniformly continuous family of uniformly entire functions.

We first prove the uniform continuous part of $f(k,T)$. Since

$$\prod_{i \geq j} f_{k,i}(\pi^i T) \equiv 1 \pmod{\pi^j},$$

we deduce that

$$f(k,T) \equiv f_j(k,T) \pmod{\pi^j}.$$

This implies that

(5.2)        $$\|f(k_1,T) - f(k_2,T)\| \leq \max(\|f_j(k_1,T) - f_j(k_2,T)\|, p^{-j}).$$

Since the family $f_j(k,T)$ parametrized by $k$ is uniformly continuous, there is a positive number $\delta(j)$ such that whenever

(5.3)                                $$\|k_1 - k_2\| < \delta(j),$$

we have the inequality

$$\|f_j(k_1,T) - f_j(k_2,T)\| < p^{-j}.$$

By (5.2), we conclude that whenever $k_1$ and $k_2$ satisfy (5.3), we have the bound

$$\|f(k_1,T) - f(k_2,T)\| \leq p^{-j}.$$

As the number $p^{-j}$ can be arbitrarily small, this proves that the family $f(k,T)$ is uniformly continuous.

Next, we prove the uniform entire part of $f(k,T)$. Write

$$f(k,T) = \sum_{m \geq 0} f_m(k) T^m,$$

$$f_j(k,T) = \sum_{m \geq 0} f_{m,j}(k) T^m,$$

$$\prod_{i \geq j} f_{k,i}(\pi^i T) = \sum_{m \geq 0} a_{m,j}(k) \pi^{jm} T^m,$$

where $a_{m,j}(k) \in R$. From the product formula

$$f(k,T) = f_j(k,T) \prod_{i \geq j} f_{k,i}(\pi^i T),$$

we deduce the relations

(5.4)        $$f_m(k) = \sum_{r=0}^{m} f_{r,j}(k) a_{m-r,j}(k) \pi^{j(m-r)}, \; m = 0, 1, 2, \cdots.$$

For any given positive constant $C$, we claim that

$$\lim_{m \to \infty} \inf \frac{\inf_k \operatorname{ord}_\pi f_m(k)}{m} \geq C.$$

In fact, take and fix a positive integer $j$ such that $j \geq 2C$. Since the family $f_j(k, T)$ parametrized by $k$ is uniformly entire, there is a positive integer $m_1$ such that for all $m > m_1$, we have the inequality

$$(5.5) \qquad \operatorname{ord}_\pi f_{m,j}(k) \geq 2Cm,$$

uniformly for all $k$. By (5.4), we deduce that for all $m > 2m_1$,

$$(5.6) \qquad \operatorname{ord}_\pi f_m(k) \geq \min_{0 \leq r \leq m} (\operatorname{ord}_\pi f_{r,j}(k) + j(m - r)).$$

Assume now that $m > 2m_1$. If $r \leq m/2$, then

$$(5.7) \qquad j(m - r) \geq j\frac{m}{2} \geq Cm.$$

If $r \geq m/2 \geq m_1$, then, by (5.5), we have

$$(5.8) \qquad \operatorname{ord}_\pi f_{r,j}(k) \geq 2Cr \geq Cm.$$

The inequalities in (5.6)-(5.8) show that the claim is true. Since $C$ can be taken to be arbitrarily large, we deduce that

$$\lim_{m \to \infty} \inf \frac{\inf_k \operatorname{ord}_\pi f_m(k)}{m} = \infty.$$

This means that the family $f(k, T)$ parametrized by $k$ is uniformly entire. The proof is complete. $\qquad \square$

## 6. Hodge-Newton decomposition

We now turn to describing Dwork's conjecture in the ordinary case. First, we discuss in more explicit detail the nuclear condition of the nuclear map $\phi$ with respect to a formal basis $\vec{e}$ of $M$. We shall introduce several notions attached to a given formal basis.

**Definition 6.1.** Let $(M, \phi)$ be a nuclear $\sigma$-module with a formal basis $\vec{e}$. For each integer $1 \leq i < \infty$, let $d_i$ be the smallest positive integer $d$ such that for all $j > d$, we have

$$\phi(e_j) \equiv 0 \pmod{\pi^i}.$$

This is a finite integer for each $i$ since

$$\lim_j \|\phi(e_j)\| = 0.$$

For $0 \leq i < \infty$, define

$$h_i = d_{i+1} - d_i, \ d_0 = 0.$$

The sequence $h = h(\vec{e}) = \{h_0, h_1, \cdots\}$ is called the **basis sequence** of $\phi$ with respect to the formal basis $\vec{e}$. Let $M_{(i)}$ be the Banach $A_0$-module with the formal basis $\{e_{d_i+1}, e_{d_i+2}, \cdots\}$. Then we have a decreasing **basis filtration** of Banach $A_0$-submodules:

$$M = M_{(0)} \supset M_{(1)} \supset \cdots \supset M_{(j)} \supset \cdots \supset 0, \ \cap_j M_{(j)} = 0,$$

where each $M_{(i)}$ has a formal basis, each quotient $M_{(i)}/M_{(i+1)}$ is a finite free $A_0$-module and

$$\phi(M_{(i)}) \subseteq \pi^i M, \ \operatorname{rank}(M_{(i)}/M_{(i+1)}) = h_i.$$

Equivalently, the matrix $B$ of $\phi$ with respect to the formal basis $\vec{e}$ (defined by $\phi(\vec{e}) = \vec{e}B$) is of the form

$$B = (B_0, \pi B_1, \pi^2 B_2, \cdots, \pi^i B_i, \cdots),$$

where each block matrix $B_i$ is a matrix over $A_0$ with $h_i$ (finitely many) columns and the last column of the $B_i$ is not divisible by $\pi$.

Next, we introduce the basis polygon and Newton polygon of a nuclear $\sigma$-module $(M, \phi)$ with a formal basis $\vec{e}$. Denote by $\text{ord}(\phi)$ to be the greatest integer $i$ such that

$$\phi \equiv 0 \ (\text{mod} \pi^i).$$

We say that $(M, \phi)$ is divisible by $\pi$ if $\phi$ is divisible by $\pi$. If $\phi$ is divisible by $\pi$, then we can write $\phi = \pi \phi_1$. In this case, one checks that

$$L(\phi, T) = L(\phi_1, \pi T).$$

Thus, for our purpose of studying L-functions, we may assume that $(M, \phi)$ is not divisible by $\pi$.

**Definition 6.2.** Let $(M, \phi)$ be a nuclear $\sigma$-module with a formal basis $\vec{e}$. Let $h = h(\vec{e}) = (h_0, h_1, \cdots)$ be the basis sequence of $\phi$ with respect to the formal basis $\vec{e}$. We define the **basis polygon** $P(\vec{e})$ of $(M, \phi)$ with respect to $\vec{e}$ to be the convex closure in the plane of the following lattice points:

$$(0, 0), (d_0, 0), (d_1, h_1), (d_2, h_1 + 2h_2), \cdots, (d_i, h_1 + 2h_2 + \cdots + ih_i), \cdots.$$

Namely, the basis polygon $P(\vec{e})$ is the polygon with a side of slope $i$ and horizontal length $h_i$ for every integer $0 \le i < \infty$.

Note that our definition of the basis polygon $P(\vec{e})$ and the basis sequence $h(\vec{e})$ depends on the given formal basis $\vec{e}$. The basis polygon and the basis sequence are somewhat similar to the Hodge polygon and the Hodge numbers. But they are different in general, even in the finite rank case. The same $\sigma$-module with different formal bases will give rise to different basis polygons and different basis sequences.

**Definition 6.3.** For each closed point $\bar{x} \in \mathbf{A}^n$ over $\mathbf{F}_q$, the Newton polygon of the nuclear $\sigma$-module $(M, \phi)$ at $\bar{x}$ is the Newton polygon of the entire characteristic series $\det(I - \phi_x^{\deg(\bar{x})} T)$ defined with respect to the valuation $\text{ord}_{\pi^{\deg(\bar{x})}}$.

A standard argument shows that the Newton polygon at each $\bar{x}$ lies on or above any basis polygon $P(\vec{e})$. This is known to be true in the finite rank case. In the infinite rank case, the proof is similar.

**Definition 6.4.** Let $\vec{e}$ be a formal basis of a nuclear $\sigma$-module $(M, \phi)$. The basis $\vec{e}$ is called **ordinary** if the basis polygon $P(\vec{e})$ of $(M, \phi)$ coincides with the Newton polygon of each fibre $(M, \phi)_x$, where $x$ is the Teichmüller lifting of the closed point $\bar{x} \in \mathbf{A}^n/\mathbf{F}_q$. The basis filtration attached to $\vec{e}$ is called ordinary if $\vec{e}$ is ordinary. Similarly, the basis $\vec{e}$ is called ordinary up to slope $j$ side if the basis polygon $P(\vec{e})$ coincides with the Newton polygon of each fibre for all sides up to slope $j$. In particular, $\vec{e}$ is said to be ordinary at the slope zero side if the horizontal side of the Newton polygon at each fibre is of length $h_0(\vec{e})$. We say that $(M, \phi)$ is ordinary (resp. ordinary up to slope $j$ side, resp. ordinary at the slope zero side) if it has a formal basis which is ordinary (resp. ordinary up to slope $j$ side, resp. ordinary at the slope zero side).

An important property is that the category of ordinary nuclear $\sigma$-modules is closed under direct sum, formal tensor product, formal symmetric power and formal exterior power. This property is not hard to prove. It follows from the proof of the Hodge-Newton decomposition. Alternatively, one could check directly. For example, let $(M_1, \phi_1)$ be a nuclear $\sigma$-module with a formal basis $\vec{e}$ and with its associated basis sequence $h = (h_0, h_1, \cdots)$. Let $(M_2, \phi_2)$ be a nuclear $\sigma$-module with a formal basis $\vec{f}$ and with its associated basis sequence $g = (g_0, g_1, \cdots)$. Then the direct sum $(M_1 \oplus M_2, \phi_1 \oplus \phi_2)$ is a nuclear $\sigma$-module with the formal basis $\vec{e} \oplus \vec{f}$. Its associated basis sequence is given by

$$h \oplus g = (h_0 + g_0, h_1 + g_1, \cdots).$$

The formal tensor product $(M_1 \otimes M_2, \phi_1 \otimes \phi_2)$ is a nuclear $\sigma$-module with the formal basis $\vec{e} \otimes \vec{f}$. Its associated basis sequence is given by

$$h \otimes g = (h_0 g_0, h_0 g_1 + h_1 g_0, h_0 g_2 + h_1 g_1 + h_2 g_0, \cdots).$$

The formal symmetric square $(\mathrm{Sym}^2 M_1, \mathrm{Sym}^2 \phi_1)$ is a nuclear $\sigma$-module with the formal basis $\mathrm{Sym}^2 \vec{e}$. Its associated basis sequence is given by

$$\mathrm{Sym}^2 h = (\frac{h_0^2 + h_0}{2}, h_0 h_1, h_0 h_2 + \frac{h_1^2 + h_1}{2}, h_0 h_3 + h_1 h_2, h_0 h_4 + h_1 h_3 + \frac{h_2^2 + h_2}{2}, \cdots).$$

The formal exterior square $(\wedge^2 M_1, \wedge^2 \phi_1)$ is a nuclear $\sigma$-module with the formal basis $\wedge^2 \vec{e}$. Its associated basis sequence is given by

$$\wedge^2 h = (\frac{h_0^2 - h_0}{2}, h_0 h_1, h_0 h_2 + \frac{h_1^2 - h_1}{2}, h_0 h_3 + h_1 h_2, h_0 h_4 + h_1 h_3 + \frac{h_2^2 - h_2}{2}, \cdots).$$

If $\vec{e}$ is an ordinary basis of $(M_1, \phi_1)$ and if $\vec{f}$ is an ordinary basis of $(M_2, \phi_2)$, one checks that $\vec{e} \oplus \vec{f}, \vec{e} \otimes \vec{f}, \mathrm{Sym}^2 \vec{e}$ and $\wedge^2 \vec{e}$ is, respectively, an ordinary basis of $(M_1 \oplus M_2, \phi_1 \oplus \phi_2)$, $(M_1 \otimes M_2, \phi_1 \otimes \phi_2)$, $(\mathrm{Sym}^2 M_1, \mathrm{Sym}^2 \phi_1)$ and $(\wedge^2 M_1, \wedge^2 \phi_1)$. The associated basis sequence is respectively $h \oplus g$, $h \otimes g$, $\mathrm{Sym}^2 h$ and $\wedge^2 h$. A systematic study of the Newton polygons and the Hodge polygons in finite rank case is given in [9].

Since $A_0$ is complete, any contraction mapping from $A_0$ into itself has a fixed point. This fact is all one needs to prove the following extension of the Hodge-Newton decomposition, first proved by Dwork [5] in the finite rank case. The infinite rank case was already considered and used in Theorem 2.4 of [12] to study the Newton polygon of Fredholm determinants arising from L-functions of exponential sums.

**Lemma 6.5.** *Let $(M, \phi)$ be a nuclear $\sigma$-module ordinary at the slope zero side. Then there is a finite free $A_0$-submodule $M_0$ of rank $h_0$ transversal to the ordinary basis filtration:*

$$M = M_0 \oplus M_{(1)},$$

*such that $\phi$ is stable on $M_0$ and $(M_0, \phi)$ is a unit root $\sigma$-module of rank $h_0$.*

In the situation of Lemma 6.5, we shall denote the sub $\sigma$-module $(M_0, \phi)$ by $(U, \phi_{\mathrm{unit}})$. We say that the unit root $\sigma$-module $(U, \phi_{\mathrm{unit}})$ is embedded in the ambient ordinary $\sigma$-module $(M, \phi)$. Alternatively, we shall say that $(U, \phi_{\mathrm{unit}})$ is the unit root part (or the slope zero part) of $(M, \phi)$. The generalized form of Dwork's unit root conjecture in our current setting is then the following

**Conjecture 6.6.** *Let $(M, \phi)$ be an overconvergent nuclear $\sigma$-module ordinary at the slope zero side. Let $(U, \phi_{\text{unit}})$ be the unit root part of $(M, \phi)$. Then for every integer $k$, the unit root zeta function $L(\phi_{\text{unit}}^k, T)$ is $p$-adic meromorphic.*

Even though the ambient $\sigma$-module $(M, \phi)$ is assumed to be overconvergent, its unit root part $\phi_{\text{unit}}$ (obtained by solving the fixed point of a contraction map) will no longer be overconvergent in general. Our main result of this paper says that Conjecture 6.6 is true if $\phi_{\text{unit}}$ has rank one. More precisely, we have

**Theorem 6.7.** *Let $(M, \phi)$ be an **overconvergent** nuclear $\sigma$-module ordinary at the slope zero side. Let $(U, \phi_{\text{unit}})$ be the unit root part of $(M, \phi)$. Assume that $\phi_{\text{unit}}$ has rank one. Let $\varphi$ be another overconvergent nuclear $\sigma$-module. Then the family of $L$-functions $L(\phi_{\text{unit}}^k \otimes \varphi, T)$ parametrized by integers $k$ in any given residue class modulo $(q-1)$ is a strong family of $p$-adic meromorphic functions.*

This result together with the result in [17] shows that Conjecture 6.6 is also true in higher rank case if $\phi$ is either of finite rank, or ordinary at every slope. However, the general form of Conjecture 6.6 for non-ordinary infinite rank $\phi$ is not proved yet. One needs to extend the method of [17] to infinite rank setup and handle some additional difficulties in the non-ordinary case for which we hope to overcome in a later paper.

Returning to Theorem 6.7. Let now $(M, \phi)$ be an overconvergent nuclear $\sigma$-module ordinary at the slope zero side. Assume that the unit root part $\phi_{\text{unit}}$ of $\phi$ has rank one. Let $\vec{e}$ be a formal overconvergent row basis of $M$ over $A_0$ which is ordinary at the slope zero side. This means that the matrix $B$ of $\phi$ with respect to $\vec{e}$ is overconvergent and the matrix $B$ has the shape

$$(6.1) \qquad B = \begin{pmatrix} B_{00} & \pi B_{01} \\ B_{10} & \pi B_{11} \end{pmatrix},$$

where $B_{00}$ is an invertible element of $A$, and the other $B_{ij}$ are matrices over $A$. Note that by our convention, the whole matrix $B$ (not just each of its entries) is also overconvergent.

**Definition 6.8.** We shall say that $\phi$ is in **normalized form** with respect to the formal overconvergent basis $\vec{e}$, if for all $i > 1$, we have

$$\phi(e_i) \equiv 0 (\text{mod } \pi),$$

and for $i = 1$ we have

$$\phi(e_1) \equiv e_1 (\text{mod } \pi).$$

In terms of the overconvergent matrix $B$ in (6.1), this means that $B_{00}$ is a 1-unit in $A$ and $B_{10}$ is divisible by $\pi$.

Since our base space is $\mathbf{A}^n$, under the assumption of Theorem 6.7, we can easily write

$$\phi = \zeta \otimes \eta,$$

where $\zeta$ is a constant (can be taken to be a root of unity) rank one $\sigma$-module and $\eta$ is a normalized overconvergent nuclear $\sigma$-module. This decomposition shows that in proving Theorem 6.7, we can assume that $\phi$ is already in normalized form. Theorem 6.7 will be proved in next section.

## 7. Limiting $\sigma$-modules

We now turn to constructing the limiting $\sigma$-module which will give us an explicit formula for the unit root L-function $L(\phi_{\mathrm{unit}}, T)$ in terms of overconvergent nuclear $\sigma$-modules. Theorem 6.7 then follows immediately.

Let $\phi_{\mathrm{unit}}$ be the slope zero part of $\phi$. A simple limiting argument applying to each Euler factor shows that

$$L(\phi_{\mathrm{unit}}^k, T) = \lim_{m \to \infty} L(\phi^{k+(q^{m!}-1)p^m}, T).$$

Combining with the decomposition formula in Lemma 2.12, we obtain a limiting formula for the unit root L-function

$$L(\phi_{\mathrm{unit}}^k, T) = \lim_{m \to \infty} \prod_{i \geq 1} L(\mathrm{Sym}^{k+(q^{m!}-1)p^m-i}\phi \otimes \wedge^i \phi, T)^{(-1)^{i-1}i}.$$

At this point, of course, we do not know if each L-factor of the right side has a limit when $m \to \infty$. Even if it has a limit, we do not know if the limit would be a meromorphic function. The purpose of this section is to answer all these questions in the affirmative in the case that $\phi_{\mathrm{unit}}$ is of rank one.

From now on in this section, our standing assumption is that $(M, \phi)$ is a nuclear overconvergent $\sigma$-module which is normalized as in Definition 6.8. Since $\phi_{\mathrm{unit}}$ is of rank one, we can replace the sequence $k + (q^{m!} - 1)p^m$ by a notationally simpler sequence such as the sequence $k_m = k + p^m$. The above limiting formula then becomes

**Lemma 7.1.** *Let $\phi$ be normalized as in Definition 6.8. Then,*

$$L(\phi_{\mathrm{unit}}^k, T) = \lim_{m \to \infty} \prod_{i \geq 1} L(\mathrm{Sym}^{k+p^m-i}\phi \otimes \wedge^i \phi, T)^{(-1)^{i-1}i}.$$

Thus, the key is to understand the limit function of the sequence of L-functions $L(\mathrm{Sym}^{k+p^m}\phi, T)$ as $m$ varies. In this section, we show that there is an overconvergent nuclear $\sigma$-module $(M_{\infty,k}, \phi_{\infty,k})$ called the limiting $\sigma$-module of the sequence $(\mathrm{Sym}^{k+p^m}M, \mathrm{Sym}^{k+p^m}\phi)$ such that

$$(7.1) \qquad \lim_{m \to \infty} L(\mathrm{Sym}^{k+p^m}\phi, T) = L(\phi_{\infty,k}, T).$$

This shows that the limit on the left side not only exists, but is also the L-function of a nuclear overconvergent $\sigma$-module and hence a meromorphic function. We shall show that the family $(M_{\infty,k}, \phi_{\infty,k})$ parametrized by $k$ is a family of uniformly overconvergent nuclear $\sigma$-modules. Furthermore, as $k$ varies $p$-adically, the family of L-functions $L(\phi_{\infty,k}, T)$ parametrized by $k$ turns out to be uniformly continuous in $k$. This gives the stronger uniform assertion stated in Theorem 6.7.

Let $k$ be a positive integer. We want to identify the $k$-th formal symmetric power $(\mathrm{Sym}^k M, \mathrm{Sym}^k \phi)$ with another nuclear overconvergent $\sigma$-module $(M_k, \phi_k)$ which is easier to take the limit as the integer $k$ varies $p$-adically. Let $\vec{e} = \{e_1, e_2, \cdots\}$ be a formal overconvergent basis of $(M, \phi)$ ordinary at the slope zero side such that its matrix satisfies the condition in (6.1). Define $M_k$ to be the Banach $A_0$-module with the formal basis $\vec{f}(k)$:

$$\{f_{i_1} f_{i_2} \cdots f_{i_r} | 2 \leq i_1 \leq i_2 \leq \cdots \leq i_r, \ 0 \leq r \leq k\},$$

where we think of 1 (corresponding to the case $r = 0$) as the first basis element of $\vec{f}(k)$. There is an isomorphism of Banach $A_0$-modules between $\mathrm{Sym}^k M$ and $M_k$.

This isomorphism is given by the map:

$$\Upsilon : e_1^{k-r} e_{i_1} \cdots e_{i_r} \longrightarrow f_{i_1} \cdots f_{i_r}, \ 0 \leq r \leq k.$$

Thus,

$$\Upsilon(e_1) = 1, \ \Upsilon(e_i) = f_i \ (i > 1).$$

We can give a nuclear overconvergent $\sigma$-module structure on $M_k$. The semi-linear map $\phi_k$ acting on $M_k$ is given by the pull back $\Upsilon \circ \mathrm{Sym}^k \phi \circ \Upsilon^{-1}$ of $\mathrm{Sym}^k \phi$ acting on $\mathrm{Sym}^k M$. Namely,

$$(7.2) \qquad \phi_k(f_{i_1} \cdots f_{i_r}) = \Upsilon(\phi(e_1^{k-r})\phi(e_{i_1}) \cdots \phi(e_{i_r})).$$

Under the identification of $\Upsilon$, the map $\mathrm{Sym}^k \phi$ becomes $\phi_k$. Thus, the $k$-th symmetric power $(\mathrm{Sym}^k M, \mathrm{Sym}^k \phi)$ is identified with $(M_k, \phi_k)$. It follows that

$$L(\mathrm{Sym}^k \phi, T) = L(\phi_k, T).$$

In order to take the limit of the sequence of modules $M_k$, we need to create a bigger space $M_\infty$ which contains each $M_k$ as a sub-module.

Let $A_0[[M]]$ be the formal symmetric algebra as constructed in section 2. The big space $M_\infty$ is defined to be the formal power series ring over $A_0$ in the variables $\{f_2, f_3, \cdots, \}$:

$$M_\infty = A_0[[f_2, f_3, \cdots]].$$

Alternatively, as the referee noted, the big ring $M_\infty$ is the quotient ring of $A_0[[M]]$ by the principal (hence closed) ideal generated by $e_1 - 1$. That is,

$$M_\infty = A_0[[M]]/(e_1 - 1).$$

This is naturally a Banach $A_0$-module with respect to the norm

$$\|C\| = \min_{c \in C} \|c\|,$$

where $C$ is a coset in $M_\infty$. The module $M_\infty$ has the countable formal basis given by $\vec{f}$:

$$(7.3) \qquad \{f_{i_1} f_{i_2} \cdots f_{i_r} | 2 \leq i_1 \leq i_2 \leq \cdots \leq i_r, \ 0 \leq r\},$$

where, again, we think of 1 (corresponding to the case $r = 0$) as the first basis element of $\vec{f}$. Note that there is no upper bound on $r$ in the definition of the basis $\vec{f}$. Each $M_k$ is a submodule of $M_\infty$. We have the inclusion

$$0 \subset M_1 \subset M_2 \subset \cdots \subset M_\infty.$$

The map $\Upsilon$ on $\mathrm{Sym}^k M$ for various $k$ fits together to give the natural reduction map

$$\Upsilon : A_0[[M]] \longrightarrow M_\infty$$

denoted by the same notation. Precisely,

$$\Upsilon : e_1^{k-r} e_{i_1} \cdots e_{i_r} \longrightarrow f_{i_1} \cdots f_{i_r}, \ 2 \leq i_1 \leq \cdots \leq i_r,$$

where we set $e_1^{k-r} = 0$ if $k < r$. Note that the $\sigma$-linear map $\mathrm{Sym}(\phi)$ acting on the formal symmetric algebra $A_0[[M]]$ is not nuclear. For each positive integer $k$, we extend the map $\phi_k$ on $M_k$ to a $\sigma$-linear map on $M_\infty$ by

$$\phi_k(f_{i_1} \cdots f_{i_r}) = \rho(\mathrm{Sym}(\phi)(e_1^{k-r} e_{i_1} \cdots e_{i_r})),$$

where we set $e_1^{k-r} = 0$ if $k < r$. Since $(M_k, \phi_k)$ is a nuclear overconvergent $\sigma$-module and $\phi_k$ induces the zero map on the quotient module $M_\infty/M_k$, it is easy to check that $(M_\infty, \phi_k)$ is also nuclear and overconvergent for every positive integer $k$. The L-function is the same whether we view $\phi_k$ as acting on the submodule $M_k$ or on the big module $M_\infty$. We shall view $\phi_k$ as acting on $M_\infty$.

To get the limiting $\sigma$-module, we need to define another nuclear $\sigma$-linear map $\phi_{\infty,k}$ on $M_\infty$ for each integer $k$. For this purpose, we take a sequence of positive integers $k_m$ (for instance $k_m = k + p^m$) such that $\lim_m k_m = \infty$ as integers and $\lim_m k_m = k$ as $p$-adic integers. Define $\phi_{\infty,k}$ by the following limiting formula:

$$
\begin{aligned}
\phi_{\infty,k}(f_{i_1} \cdots f_{i_r}) &= \lim_{m\to\infty} \Upsilon(\mathrm{Sym}(\phi)(e_1^{k_m-r} e_{i_1} \cdots e_{i_r})) \\
&= \big( \lim_{m\to\infty} \Upsilon(\phi(e_1)^{k_m-r})\big) \Upsilon(\phi(e_{i_1}) \cdots \phi(e_{i_r})).
\end{aligned}
$$

To show that this is well defined, we need to show that the above first limiting factor exists. Write

$$(7.4) \qquad \phi(e_1) = e_1 + \pi e,$$

where $e$ is an element of $M$. By the binomial theorem, for each positive integer $j$, we have

$$
\begin{aligned}
\phi(e_1^j) &= (e_1 + \pi e)^j \\
&= e_1^j + \binom{j}{1}\pi e_1^{j-1}e + \binom{j}{2}\pi^2 e_1^{j-2}e^2 + \cdots.
\end{aligned}
$$

This identity implies that

$$(7.5) \qquad \lim_{m\to\infty} \Upsilon(\phi(e_1)^{k_m-r}) = 1 + \binom{k-r}{1}\pi\Upsilon(e) + \binom{k-r}{2}\pi^2\Upsilon(e)^2 + \cdots.$$

Thus,

$$(7.6) \qquad \phi_{\infty,k}(f_{i_1} \cdots f_{i_r}) = (1 + \pi\Upsilon(e))^{k-r}\Upsilon(\phi(e_{i_1}) \cdots \phi(e_{i_r})).$$

It follows that the map $\phi_{\infty,k}$ is well defined for every integer $k$ and it is independent of the choice of our chosen sequence $k_m$. Furthermore, using (7.6), we see that $\phi_{\infty,k}$ makes sense for any $p$-adic integer $k$, not necessarily usual positive integers. The map $\phi_{\infty,k}$ is a $\sigma$-linear ring endomorphism of $M_\infty$. Our construction of the limiting $\sigma$-modules $(M_\infty, \phi_{\infty,k})$ depends on the choice of the basis $\vec{e}$. As Coleman has observed, the limiting $\sigma$-module can in fact be constructed in a functorial way in some sense.

The relation between $\phi_{\infty,k}$ and $\phi_k$ will be discussed later. First, we try to understand the limiting $\sigma$-module $\phi_{\infty,k}$. We have the following basic result.

**Theorem 7.2.** *Let $(M, \phi)$ be a nuclear overconvergent $\sigma$-module which is ordinary at the slope zero side. Assume that $\phi_{\mathrm{unit}}$ has rank one and $\phi$ is in normalized situation as in (6.1). Then, the family $(M_\infty, \phi_{\infty,k})$ parametrized by $p$-adic integer $k$ is a uniformly continuous family of uniformly overconvergent nuclear $\sigma$-modules. In particular, by Corollary 5.9, the family of L-functions $L(\phi_{\infty,k}, T)$ parametrized by $p$-adic integer $k \in \mathbf{Z}_p$ is a strong family of meromorphic functions.*

*Proof.* First, we check that the family of L-functions $L(\phi_{\infty,k}, T)$ is uniformly continuous in $k$ with respect to the $p$-adic topology. By our definition of $\phi_{\infty,k}$ in (7.6), it suffices to show that

$$\|(1 + \pi\Upsilon(e))^{k_1-r} - (1 + \pi\Upsilon(e))^{k_2-r}\| \le c\|k_1 - k_2\|$$

uniformly for all integers $k_1$ and $k_2$, where $c$ is some positive constant depending only on $\pi$. But this follows from the binomial theorem. The uniform continuity is established. The following two lemmas show that the family $\phi_{\infty,k}$ is uniformly nuclear and uniformly overconvergent. The theorem is then proved. $\blacksquare$

**Lemma 7.3.** *The family $\phi_{\infty,k}$ is uniformly nuclear. That is, for any given positive integer $C$, there is another positive integer $C_1 > 0$ such that the inequality*

$$\mathrm{ord}_\pi \phi_{\infty,k}(f_{i_1} \cdots f_{i_r}) \geq C$$

*holds uniformly for all $k$ and all $2 \leq i_1 \leq i_2 \leq \cdots \leq i_r$ with $i_1 + \cdots + i_r > C_1$.*

*Proof.* Because $\phi$ is nuclear, there is a finite integer $d$ depending on the given positive integer $C$ such that

$$\mathrm{ord}_\pi \phi(e_i) \geq C$$

holds for all $i > d$. Thus, by equation (7.6), we may restrict our attention of the indices $\{i_1, \cdots, i_r\}$ to the range

$$2 \leq i_1 \leq i_2 \leq \cdots \leq i_r \leq d.$$

Since $\phi$ is normalized with respect to $\vec{e}$, $\phi(e_i)$ is divisible by $\pi$ for every $i \geq 2$. This shows that the product $\phi(e_{i_1})\phi(e_{i_2}) \cdots \phi(e_{i_r})$ is divisible by $\pi^C$ if $r \geq C$. By (7.6) again, we may assume that $r \leq C$. Now, there are only finitely many integer tuples $\{i_1, \cdots, i_r\}$ satisfying the conditions

$$0 \leq r \leq C, \ 2 \leq i_1 \leq i_2 \leq \cdots \leq i_r \leq d.$$

It is clear that if

$$i_1 + \cdots + i_r > dC,$$

then either $r > C$ or $i_r > d$. In either case, we have the inequality

$$\mathrm{ord}_\pi \phi_{\infty,k}(f_{i_1} \cdots f_{i_r}) \geq C$$

uniformly for all $k$. The Lemma holds with $C_1 = dC$. $\blacksquare$

To finish the proof of Theorem 7.2, we need to prove

**Lemma 7.4.** *The family $(M_\infty, \phi_{\infty,k})$ is uniformly overconvergent.*

*Proof.* For rational numbers $b > 0$ and $c$, we define

$$M_\infty(b,c) = \{ \sum_{2 \leq i_1 \leq i_2 \leq \cdots \leq i_r} a_{i_1, \cdots, i_r} f_{i_1} \cdots f_{i_r} | a_{i_1, \cdots, i_r} \in L(b,c) \}.$$

This is a complete submodule of the ring $M_\infty$. Since $L(b,0)$ is a ring, it follows that $M(b,0)$ is also a ring containing $L(b,0)$ as a subring. In fact, $M_\infty(b,0)$ is an infinite dimensional commutative $L(b,0)$-algebra. More generally,

$$M_\infty(b,c_1)M_\infty(b,c_2) \subset M_\infty(b,c_1+c_2).$$

Now, since $\phi$ is overconvergent with respect to the basis $\vec{e}$, there are rational numbers $b > 0$ and $c$ such that for all $j \geq 1$,

$$\phi(e_j) \in M(b,c).$$

Since $\phi(e_j)$ is divisible by $\pi$ for each $j \geq 2$, a simple geometric argument in the plane $\mathbf{R}^2$ shows that we can choose $b$ sufficiently small that

$$\phi(e_j) \in M(b,0), \ j \geq 2.$$

This implies that

(7.7)                             $\Upsilon(\phi(e_j)) \in M_\infty(b,0), \ j \geq 2.$

For $j = 1$, we have

$$\phi(e_1) = e_1 + \pi e,$$

for some $e \in M$. We can choose our $b$ smaller if necessary so that

$$\pi e \in M(b, 0).$$

This implies that

$$\Upsilon(\pi e) \in M_\infty(b, 0).$$

Since $M_\infty(b, 0)$ is a ring, we deduce that for all positive integers $j$,

$$\Upsilon(\pi e)^j \in M_\infty(b, 0).$$

The binomial theorem shows that for any $p$-adic integer $k \in \mathbf{Z}_p$, we have

$$(7.8) \qquad (1 + \Upsilon(\pi e))^k = 1 + \binom{k}{1}\Upsilon(\pi e) + \binom{k}{2}\Upsilon(\pi e)^2 + \cdots \in M_\infty(b, 0).$$

By (7.6)-(7.8), we conclude that for all $2 \le i_1 \le i_2 \le \cdots \le i_r$, we always have

$$\phi_{\infty, k}(f_{i_1} \cdots f_{i_r}) \in M_\infty(b, 0),$$

uniformly for all $p$-adic integers $k$. The Lemma is proved. $\qquad\blacksquare$

We have finished the proof of Theorem 7.2. It remains to connect the L-function $L(\phi_{\infty, k}, T)$ in Theorem 7.2 with the desired unit root L-function $L(\phi_{\text{unit}}^k, T)$. This is done by a simple limiting argument.

**Lemma 7.5.** *Let $k$ be an integer and let $k_m = k + p^m > 0$. Then as endomorphisms on $M_\infty$, we have the congruence*

$$\phi_{k_m} \equiv \phi_{\infty, k} \pmod{\pi^{\min(k_m, m)}}.$$

*Proof.* We need to show that for all indices $2 \le i_1 \le i_2 \le \cdots \le i_r$,

$$\phi_{k_m}(f_{i_1} \cdots f_{i_r}) \equiv \phi_{\infty, k}(f_{i_1} \cdots f_{i_r}) \pmod{\pi^{\min(k_m, m)}}.$$

If $r > k_m$, the left side is zero and the right side satisfies

$$\begin{aligned} \mathrm{ord}_\pi(\phi_{\infty, k}(f_{i_1} \cdots f_{i_r})) &\ge& \mathrm{ord}_\pi(\phi(e_{i_1}) + \cdots + \mathrm{ord}_\pi(\phi(e_{i_r})) \\ &\ge& r \ge k_m. \end{aligned}$$

The congruence is indeed true. Assume now that $r \le k_m$. By (7.6), it suffices to check that

$$\Upsilon(\phi(e_1^{k_m - r})) \equiv \lim_\ell \Upsilon(\phi(e_1^{k + p^\ell - r})) \pmod{\pi^m}.$$

Equivalently, one needs to check

$$(7.9) \qquad (1 + \pi\Upsilon(e))^{k - r + p^m} \equiv (1 + \pi\Upsilon(e))^{k - r} \pmod{\pi^m}.$$

But (7.9) is a consequence of the binomial theorem. The Lemma is proved. $\qquad\blacksquare$

**Corollary 7.6.** *Let $k$ be an integer and let $k_m = k + p^m$. Then, we have the limiting formula*

$$\begin{aligned} \lim_{m \to \infty} L(\mathrm{Sym}^{k_m}\phi, T) &=& \lim_{m \to \infty} L(\phi_{k_m}, T) \\ &=& L(\phi_{\infty, k}, T). \end{aligned}$$

*More generally, for each integer $i \geq 0$, we have the limiting formula*

$$
\begin{aligned}
\lim_{m \to \infty} L(\mathrm{Sym}^{k_m} \phi \otimes \wedge^i \phi, T) &= \lim_{m \to \infty} L(\phi_{k_m} \otimes \wedge^i \phi, T) \\
&= L(\phi_{\infty,k} \otimes \wedge^i \phi, T).
\end{aligned}
$$

Using these results, we obtain the following explicit formula for the unit root L-function.

**Theorem 7.7.** *Let $(M, \phi)$ be an **overconvergent** nuclear $\sigma$-module, ordinary at the slope zero side. Let $(U, \phi_{\mathrm{unit}})$ be the unit root part of $(M, \phi)$. Assume that $\phi_{\mathrm{unit}}$ is of rank one and in normalized situation. Denote by $\phi_{\infty,k}$ the limiting nuclear $\sigma$-module of the sequence $\mathrm{Sym}^{k+p^m} \phi$. Then for all integers $k$, we have the following explicit formula*

$$
(7.10) \qquad L(\phi_{\mathrm{unit}}^k, T) = \prod_{i \geq 1} L(\phi_{\infty,k-i} \otimes \wedge^i \phi, T)^{(-1)^{i-1}i}.
$$

*In particular, the family $L(\phi_{\mathrm{unit}}^k, T)$ of unit root L-functions parametrized by $p$-adic integer $k$ is a strong family of meromorphic functions with respect to the $p$-adic topology of $k$.*

It should be noted that the product in (7.10) is a finite product if $\phi$ is of finite rank, since we have $\wedge^i \phi = 0$ for $i$ greater than the rank of $\phi$. If $\phi$ is of infinite rank, then the product in (7.10) is an infinite product. In this case, the strong family assertion about $L(\phi_{\mathrm{unit}}^k, T)$ follows from Lemma 5.10. This is because $\phi$ is nuclear and thus $\wedge^i \phi$ is more and more divisible by $\pi$ as $i$ grows. We can then apply Lemma 5.10 to the even indexed product and the odd indexed product in (7.10), respectively.

To get the full Theorem 6.7, we need to twist the unit root family $\phi_{\mathrm{unit}}^k$ by a fixed nuclear overconvergent $\sigma$-module $\varphi$. We state this generalization here. The proof is the same as the proof of Theorem 7.7. One simply twists the whole proof by the harmless overconvergent $\varphi$. The twisted basic decomposition formula is

$$
L(\phi^k \otimes \varphi, T) = \prod_{i=1}^{\infty} L(\mathrm{Sym}^{k-i} \phi \otimes \wedge^i \phi \otimes \varphi, T)^{(-1)^{i-1}i}.
$$

The required uniform overconvergent nuclear result is already given in Theorem 7.2.

**Theorem 7.8.** *Let $(M, \psi)$ and $(N, \varphi)$ be two **overconvergent** nuclear $\sigma$-modules. Assume that the first one $\psi$ is ordinary at the slope zero side whose unit root part $\psi_{\mathrm{unit}}$ is of rank one. Write $\psi = a \otimes \phi$, where $a$ is a $p$-adic unit in $R$, and $\phi$ is in normalized situation as in (6.1). Denote by $\phi_{\infty,k}$ the limiting nuclear $\sigma$-module of the sequence $\mathrm{Sym}^{k+p^m} \phi$. Then for all integers $k$, we have the following explicit formula*

$$
L(\psi_{\mathrm{unit}}^k \otimes \varphi, T) = \prod_{i \geq 1} L(a^k \otimes \phi_{\infty,k-i} \otimes \wedge^i \phi \otimes \varphi, T)^{(-1)^{i-1}i}.
$$

*In particular, the family $L(\psi_{\mathrm{unit}}^k \otimes \varphi, T)$ of twisted unit root L-functions parametrized by integers $k$ in any residue class modulo $(q-1)$ is a strong family of meromorphic functions with respect to the $p$-adic topology of $k$.*

This completes our proof of Theorem 6.7.

## 8. Distribution of zeros and poles

First, we assume that $\phi$ is in normalized situation as in Theorem 7.7. We want to get some explicit information about the zeros and poles of the unit root L-function $L(\phi_{\text{unit}}^k, T)$. For this purpose, we need to get a non-trivial uniform lower bound for the Newton polygon of the numerator and denominator of $L(\phi_{\text{unit}}^k, T)$. By Theorem 7.7, it suffices to get a non-trivial uniform lower bound for the L-function $L(\phi_{\infty,k} \otimes \varphi, T)$, where $\varphi$ is any fixed overconvergent nuclear $\sigma$-module. This can indeed be done, since our proof of Theorem 7.7 is completely explicit.

To get cleaner estimates, we shall assume that both $\phi$ and $\varphi$ are of finite rank. This is enough for most of the immediate applications we have in mind. Thus, we shall assume that $\phi$ is of some finite rank $r \geq 1$ and in normalized situation. One could further assume that $r \geq 2$ as the case $r = 1$ reduces to the classical overconvergent case. We shall however include the trivial case $r = 1$ as well, since all of our results are also true for the exceptional case $r = 1$. This provides a useful comparison of results between the old overconvergent case and the new non-overconvergent case.

We first consider the L-function $L(\phi_{\infty,k}, T)$. This is the essential case. By the construction of the limiting $\sigma$-module in section 7, we know that the underlying module $M_\infty$ of $\phi_{\infty,k}$ is the formal power series ring over $A_0$ in the $(r-1)$-variables $\{f_2, \cdots, f_r\}$:

$$
\begin{aligned}
M_\infty &= A_0[[f_2, \cdots, f_r]] \\
&= \{ \sum_{i_2, \cdots, i_r} a_{i_2, \cdots, i_r} f_2^{i_2} \cdots f_r^{i_r} | a_{i_2, \cdots, i_r} \in A_0 \},
\end{aligned}
$$

where the $i_\ell$ run over the set of non-negative integers. The corresponding continuous dual $A_0$-module is the convergent power series ring

$$
\begin{aligned}
M_\infty^* &= A_0\{\{f_2^*, \cdots, f_r^*\}\} \\
&= \{ \sum_{i_2, \cdots, i_r} a_{i_2, \cdots, i_r} f_2^{*i_2} \cdots f_r^{*i_r} | a_{i_2, \cdots, i_r} \in A_0, \ \lim \|a_{i_2, \cdots, i_r}\| = 0 \}.
\end{aligned}
$$

Note that in the exceptional case $r = 1$, both $M_\infty$ and $M_\infty^*$ are rank one (not rank zero) $A_0$-module, corresponding to the constant term $i_2 = \cdots = i_r = 0$. In this exceptional case, one checks that the limiting $\sigma$-module is simply given by

$$
\phi_{\infty,k} = \lim_{m \to \infty} \text{Sym}^{k+p^m} \phi = \phi^k
$$

which is a rank one overconvergent $\sigma$-module.

Since $\phi$ is normalized with respect to $\vec{e}$, we have the congruence

$$
\phi(e_1) \equiv e_1 \ (\text{mod } \pi), \ \phi(e_i) \equiv 0 \ (\text{mod } \pi), \ 2 \leq i \leq r.
$$

We may thus choose a smaller $b$ and go to a totally ramified finite extension of $R$ if necessary such that

$$
\phi(e_1) \in M(b, 0), \ \phi(e_i) \in \pi M(b, 0), \ 2 \leq i \leq r.
$$

From the construction in (7.6) of the limiting $\sigma$-module, we deduce

$$
\phi_{\infty,k}(f_i) \in \pi M_\infty(b, 0), \quad 2 \leq i \leq r.
$$

Since $L(b, 0)$ is a ring, we obtain

$$
\phi_{\infty,k}(f_2^{i_2} \cdots f_r^{i_r}) \in \pi^{i_2 + \cdots + i_r} M_\infty(b, 0).
$$

Let $\Theta(k)$ be the Dwork operator on $M_\infty^*$ associated to $\phi_{\infty,k}$. Write

$$\Theta(k)(\pi^{b|v|} X^v f_2^{*j_2} \cdots f_r^{*j_r}) = \sum_{u,i_2,\cdots,i_r} G_{\{u,i\},\{v,j\}}(k) \pi^{b|u|} X^u f_2^{*i_2} \cdots f_r^{*i_r},$$

where

$$i = \{i_2, \cdots, i_r\}, \ j = \{j_2, \cdots, j_r\}.$$

Applying Lemma 5.4 to $\Theta(k)$ acting on $M_\infty^*$, we derive the uniform estimate

$$\mathrm{ord}_\pi G_{\{u,i\},\{v,j\}}(k) \geq (q-1)b|u| + |i| + c, \ |i| = i_2 + \cdots + i_r.$$

This shows that for the positive constant $c_1 = \min((q-1)b, 1)$, we have the bound

$$\mathrm{ord}_\pi G_{\{u,i\},\{v,j\}}(k) \geq c_1(|u| + |i|) + c,$$

uniformly for all $i = (i_2, \cdots, i_r), j = (j_2, \cdots, j_r)$ and $k$, where $c$ is a constant.

Consider the Newton polygon whose $\ell$-th vertex ($0 \leq \ell < \infty$) is given by

$$(8.1) \qquad \Big( \sum_{u_1,\cdots,u_n,i_2,\cdots,i_r} 1, \sum_{u_1,\cdots,u_n,i_2,\cdots,i_r} c_1(u_1 + \cdots + u_n + i_2 + \cdots + i_r) + c \Big),$$

where the $u_j$ and $i_k$ run over non-negative integers satisfying the inequality

$$u_1 + \cdots + u_n + i_2 + \cdots + i_r \leq \ell.$$

Let $P(x)$ be the real valued function on the non-negative real numbers $\mathbf{R}_{\geq 0}$ whose graph is the Newton polygon in (8.1). A standard determinant argument then shows that the Newton polygon of the Fredholm determinant $\det(I - \Theta(k)T|M_\infty^*)$ lies on or above the Newton polygon $P(x)$. Thus, we need to estimate the growth of the function $P(x)$.

Now, the first coordinate in (8.1) is just the coefficient of $t^\ell$ in the power series expansion of the function $(1-t)^{-(n+r)}$. This coefficient is, by the binomial theorem,

$$\binom{n+r+\ell-1}{n+r-1} \sim c_2(n+r+\ell-1)^{n+r-1}$$

for large $\ell$ and some positive constant $c_2$. Similarly, an easy argument shows that the second coordinate of (8.1) is bounded from below by

$$c_3(n+r+\ell-1) \sum_{u_1+\cdots+u_n+i_2+\cdots+i_n \leq \ell} 1 \geq c_4(n+r+\ell-1)^{n+r}$$

for large $\ell$ and some positive constants $\{c_3, c_4\}$. This together with the concavity of the smooth function $x^{n+r}$ implies that the Newton polygon $P(x)$ satisfies the lower bound

$$P(x^{n+r-1}) \geq c_5 x^{n+r}$$

for all large $x$ and some positive constant $c_5$. Thus, we deduce for all large $x$,

$$P(x) \geq c_5 x^{1 + \frac{1}{n+r-1}}.$$

Since $P(0) = 0$, we can choose $c_6$ to be a sufficiently large integer such that the following result holds.

**Lemma 8.1.** *There are positive constants $c_5$ and $c_6$ such that for all real numbers $x \geq 0$, we have*

$$P(x) \geq c_5 x^{1 + \frac{1}{n+r-1}} - c_6 x.$$

In this lemma, we used the linear error term $c_6 x$ for simplicity. If one wishes, one could replace the linear term $c_6 x$ in Lemma 8.1 by any positive function which grows slower than $x^{1+1/(n+r-1)}$.

The lower bound in Lemma 8.1 provides a uniform lower bound for the Newton polygon of the family of Fredholm determinants $\det(I - T\Theta(k)|M_\infty^*)$ parametrized by $k$. To handle the L-function $L(\phi_{\infty,k}, T)$, we also need to do the same for the $i$-th Dwork operator $\Theta_i(k)$ associated to $\phi_{\infty,k}$, corresponding to differential $i$-forms. Thus, let $i$ be any integer between $0$ and $n$. As before, we define the continuous dual $A_0$-module of $\Omega^i M_\infty$ by

$$
\begin{aligned}
M_{\infty,i}^* &= \mathrm{Hom}_{A_0}^{\mathrm{cont}}(\Omega^i M_\infty, \Omega^n A_0) \\
&= \{ \sum_{j_2,\cdots,j_r} a_{j_2,\cdots,j_r} f_2^{*j_2} \cdots f_r^{*j_r} | \lim_j \|a_j\| = 0 \},
\end{aligned}
$$

where

$$
a_{j_2,\cdots,j_r} \in \mathrm{Hom}_{A_0}(\Omega^i A_0, \Omega^n A_0).
$$

That is, $M_{\infty,i}^*$ consists of continuous $A_0$-linear maps from $\Omega^i M_\infty$ to $\Omega^n A_0$. Let $\Theta_i(k)$ be the $i$-th Dwork operator on $M_{\infty,i}^*$ associated to $\phi_{\infty,k}$. The same argument as above shows that we have

**Lemma 8.2.** *There are positive constants $c_5$ and $c_6$ such that the Newton polygon of $\det(I - T\Theta_i(k)|M_{\infty,i}^*)$ lies above the graph of the function*

$$
Q(x) = c_5 x^{1+\frac{1}{n+r-1}} - c_6 x
$$

*uniformly for every integer $0 \le i \le n$ and every $p$-adic integer $k$.*

So far, we only considered the essential case when $\phi$ is in normalized situation. The general case can be done in the same way, except that we have to twist the whole proof by a fixed overconvergent finite rank $\sigma$-module. The resulting bounds are the same except that we may need to use a smaller constant $c_5$ and a larger constant $c_6$. We state the result here.

**Lemma 8.3.** *Let $\varphi$ be a finite rank overconvergent $\sigma$-module. Let $\Theta_i(k,\varphi)$ be the $i$-th Dwork operator associated to the nuclear $\sigma$-module $\phi_{\infty,k} \otimes \varphi$. There are positive constants $c_5$ and $c_6$ such that the Newton polygon of $\det(I - T\Theta_i(k,\varphi))$ lies above the graph of the function*

$$
Q(x) = c_5 x^{1+\frac{1}{n+r-1}} - c_6 x
$$

*uniformly for every integer $0 \le i \le n$ and every $p$-adic integer $k$.*

Our main result in this section is

**Theorem 8.4.** *Let $(M, \psi)$ and $(N, \varphi)$ be two overconvergent finite rank $\sigma$-modules. Assume that the first one $\psi$ is of rank $r$, ordinary at the slope zero side whose unit root part $\psi_{\mathrm{unit}}$ is of rank one. Then, for any integer $k$, we can write*

$$
L(\psi_{\mathrm{unit}}^k \otimes \varphi, T) = \frac{f_1(k,T)}{f_2(k,T)},
$$

*where $f_1(k,T)$ (resp. $f_2(k,T)$) is a family of uniformly entire functions whose Newton polygons lies above the graph of the function*

$$
Q(x) = c_5 x^{1+\frac{1}{n+r-1}} - c_6 x
$$

*for some positive constants $\{c_5, c_6\}$ independent of $k$. Furthermore, if $k_1$ and $k_2$ are two integers in the same residue class modulo $(q-1)$, then we have*

$$\|f_1(k_1, T) - f_1(k_2, T)\| \le c\|k_1 - k_2\|, \ \ \|f_2(k_1, T) - f_2(k_2, T)\| \le c\|k_1 - k_2\|,$$

*for some positive constant $c$ depending only on $\pi$.*

*Proof.* Write $\psi = a \otimes \phi$, where $a$ is a $p$-adic unit in $R$, and $\phi$ is in normalized situation as in (6.1). Denote by $\phi_{\infty,k}$ the limiting nuclear $\sigma$-module of the sequence $\mathrm{Sym}^{k+p^m}\phi$. Let $\Theta_i(k-j, a^k \wedge^j \phi \otimes \varphi)$ denote the $i$-th Dwork operator associated to the nuclear $\sigma$-module $\phi_{\infty,k-j} \otimes a^k \otimes \wedge^j \phi \otimes \varphi$. Then for all integers $k$, by Theorem 7.8 and Theorem 5.8, we deduce the following explicit formula

$$
\begin{aligned}
L(\psi_{\mathrm{unit}}^k \otimes \varphi, T) &= \prod_{j=1}^r L(\phi_{\infty,k-j} \otimes a^k \otimes \wedge^j \phi \otimes \varphi, T)^{(-1)^{j-1}j} \\
&= \prod_{j=1}^r \prod_{i=0}^n \det(I - T\Theta_i(k-j, a^k \wedge^j \phi \otimes \varphi))^{(-1)^{i+j}j}.
\end{aligned}
$$

Let

$$f_1(k, T) = \prod_{j+i \text{ even}} \det(I - T\Theta_i(k-j, a^k \wedge^j \phi \otimes \varphi))^j$$

and

$$f_2(k, T) = \prod_{j+i \text{ odd}} \det(I - T\Theta_i(k-j, a^k \wedge^j \phi \otimes \varphi))^j.$$

Then, we have

$$L(\psi_{\mathrm{unit}}^k \otimes \varphi, T) = \frac{f_1(k, T)}{f_2(k, T)}.$$

The rest of the properties follows from Lemma 8.3 and Corollary 4.17. We only need to note that if $k_1$ and $k_2$ are two integers in the same residue class modulo $(q-1)$, then

$$\|\phi_{\infty,k_1} \otimes a^{k_1} \otimes \wedge^j \phi \otimes \varphi - \phi_{\infty,k_2} \otimes a^{k_2} \otimes \wedge^j \phi \otimes \varphi\| \le c\|k_1 - k_2\|$$

for some positive constant depending only on $\pi$. The theorem is proved. ∎

Now, we turn to the Gouvêa-Mazur type conjecture for the unit root L-function. Assume that we are in the situation of Theorem 8.4. By the $p$-adic Weierstrass factorization theorem, we can write

$$L(\psi_{\mathrm{unit}}^k \otimes \varphi, T) = \prod_{j \ge 1} (1 - z_j(k)T)^{\pm},$$

where $z_j(k)$ denotes a typical reciprocal zeros or poles. For a given rational non-negative rational number $s$, let

$$L_s(k, T) = \prod_{\mathrm{ord}_\pi z_j(k) = s} (1 - z_j(k)T)^{\pm}$$

denote the slope $s$ part of the above L-function. This is a rational function with coefficients in $R$. Let $d_s(k)$ (resp. $D_s(k)$) denote the degree (resp. the total degree) of the rational function $L_s(k, T)$. The function $d_s(k)$ (resp. $D_s(k)$) of two variables is called the degree (resp. the total degree) function of the meromorphic L-function $L(\psi_{\mathrm{unit}}^k \otimes \varphi, T)$. We would like to understand these two variable functions $d_s(k)$ and $D_s(k)$.

Similarly, we can define a degree function for each of the Fredholm determinant in the above formula for the L-function. For $0 \leq i \leq n$ and $1 \leq j \leq r$, we define $D_{s,i}(k-j, a^k \wedge^j \phi \otimes \varphi)$ to be the degree of the slope $s$ part of the entire Fredholm determinant $\det(I - T\Theta_i(k-j, a^k \wedge^j \phi \otimes \varphi))$. We would also like to understand these two variable functions for each fixed $i$ and $j$. They are related to $d_s(k)$ and $D_s(k)$ by the following relations:

$$d_s(k) = \sum_{j+i \text{ even}} j D_{s,i}(k-j, a^k \wedge^j \phi \otimes \varphi) - \sum_{j+i \text{ odd}} j D_{s,i}(k-j, a^k \wedge^j \phi \otimes \varphi)$$

and

$$D_s(k) \leq \sum_{j=1}^{r} \sum_{i=0}^{n} j D_{s,i}(k-j, a^k \wedge^j \phi \otimes \varphi).$$

**Theorem 8.5.** *Assume that we are in the situation of Theorem 8.4. There is a positive constant $c$ independent of $k$ such that for all real numbers $s \geq 1$, the inequality*

$$\frac{1}{s} \sum_{t \in [0,s]} D_{t,i}(k-j, a^k \wedge^j \phi \otimes \varphi) \leq c(s+1)^{n+r-2}$$

*holds uniformly for all integers $k$, all integers $0 \leq i \leq n$ and all integers $1 \leq j \leq r$.*

*Proof.* Lemma 8.3 implies that there are two positive constants $c_5$ and $c_6$ independent of $k$ such that

$$s \Big( \sum_{t \in [0,s]} D_{t,i}(k-j, a^k \wedge^j \phi \otimes \varphi) \Big)$$

$$\geq \sum_{t \in [0,s]} t D_{t,i}(k-j, a^k \wedge^j \phi \otimes \varphi)$$

$$\geq c_5 \Big( \sum_{t \in [0,s]} D_{t,i}(k-j, a^k \wedge^j \phi \otimes \varphi) \Big)^{1+\frac{1}{n+r-1}}$$

$$- c_6 \sum_{t \in [0,s]} D_{t,i}(k-j, a^k \wedge^j \phi \otimes \varphi).$$

Cancel the factor

$$\sum_{t \in [0,s]} D_{t,i}(k-j, a^k \wedge^j \phi \otimes \varphi)$$

from the above inequality (we may assume that it is non-zero, otherwise the result is trivial), we deduce

$$s \geq c_5 \Big( \sum_{t \in [0,s]} D_{t,i}(k-j, a^k \wedge^j \phi \otimes \varphi) \Big)^{\frac{1}{n+r-1}} - c_6.$$

This implies that

$$\sum_{t \in [0,s]} D_{t,i}(k-j, a^k \wedge^j \phi \otimes \varphi) \leq c(s+1)^{n+r-1}$$

uniformly. The theorem is proved.                                        $\square$

**Corollary 8.6.** *Assume that we are in the situation of Theorem 8.4. There is a positive constant $c$ independent of $k$ such that for all real numbers $s \geq 1$, the inequality*

$$\frac{1}{s} \sum_{t \in [0,s]} |d_t(k)| \leq \frac{1}{s} \sum_{t \in [0,s]} D_t(k) \leq c(s+1)^{n+r-2}$$

*holds uniformly for all $k$.*

**Corollary 8.7.** *Assume that we are in the situation of Theorem 8.4. There is a positive constant $c$ such that for all $s \geq 0$, we have the uniform polynomial bounds*

$$D_{s,i}(k-j, a^k \wedge^j \phi \otimes \varphi) \leq \sum_{t \in [0,s]} D_{t,i}(k-j, a^k \wedge^j \phi \otimes \varphi) \leq c(s+1)^{n+r-1},$$

$$\sum_{t \in [0,s]} |d_t(k)| \leq \sum_{t \in [0,s]} D_t(k) \leq c(s+1)^{n+r-1},$$

$$|d_s(k)| \leq D_s(k) \leq c(s+1)^{n+r-1}$$

*for all integers $k$, all integers $0 \leq i \leq n$ and all integers $1 \leq j \leq r$.*

The Gouvêa-Mazur type question in our situation is how the integer $d_s(k)$ varies with the parameter $k$ for a fixed slope $s$. We now turn to this question. We have the following result.

**Proposition 8.8.** *Assume that we are in the situation of Theorem 8.4. There is a positive constant $c$ such that whenever $k_1$ and $k_2$ are two integers satisfying the congruence*

$$k_1 \equiv k_2 \pmod{(q-1)p^{c[s+1]^{n+r}}},$$

*we have the equality*

$$D_{t,i}(k_1 - j, a^k \wedge^j \phi \otimes \varphi) = D_{t,i}(k_2 - j, a^k \wedge^j \phi \otimes \varphi)$$

*for all rational numbers $0 \leq t \leq s$, all integers $0 \leq i \leq n$ and all integers $1 \leq j \leq r$.*

*Proof.* First, we do have the continuity result for $\det(I - T\Theta_i(k-j, \wedge^j \phi \otimes \varphi))$. That is,

$$\| \det(I - T\Theta_i(k_1 - j, a^{k_1} \wedge^j \phi \otimes \varphi)) - \det(I - T\Theta_i(k_2 - j, a^{k_2} \wedge^j \phi \otimes \varphi)) \| \leq c|k_1 - k_2|$$

as long as $k_1$ and $k_2$ are in the same residue class modulo $(q-1)$, where $c$ is a positive constant depending only on $\pi$. Second, we have the uniform lower bound for the Newton polygon of $\det(I - T\Theta_i(k-j, a^k \wedge^j \phi \otimes \varphi))$. This uniform lower bound is given by the function

$$Q(x) = c_5 x^{1 + \frac{1}{n+r-1}} - c_6 x.$$

Taking $\nu(x) = (c_5 x^{\frac{1}{n+r-1}} - c_6)$ in the following lemma (Lemma 4.1 in [15]), we see that the proposition is proved. $\qquad \blacksquare$

**Lemma 8.9.** *Let $g_1(T)$ and $g_2(T)$ be two elements in $R[[T]]$ with $g_1(0) = g_2(0) = 1$. Let $N_i(x)$ be the function on $\mathbf{R}_{\geq 0}$ whose graph is the Newton polygon of $g_i(T)$ ($1 \leq i \leq 2$). Assume that $\nu(x)$ is a strictly increasing continuous function on $\mathbf{R}_{\geq 0}$ such that*

$$\nu(0) \leq 0, \ N_i(x) \geq x\nu(x) \ (1 \leq i \leq 2, x \geq 1), \ \lim_{x \to \infty} \nu(x) = \infty.$$

*Assume further that the function $x\nu^{-1}(x)$ is increasing on $\mathbf{R}_{\geq 0}$, where $\nu^{-1}(x)$ denotes the inverse function of $\nu(x)$ defined at least on $\mathbf{R}_{\geq 0}$ as $\nu(0) \leq 0$. For $x \geq 0$, we define the integer valued increasing function*

$$m_\nu(x) = [x\nu^{-1}(x)].$$

*If the formal power series congruence*

$$g_1(T) \equiv g_2(T) \pmod{\pi^{m_\nu(s)+1}}$$

*holds for some $s \geq 0$, then the two Newton polygons $N_i(x)$ coincide for all the sides with slopes at most $s$.*

**Corollary 8.10.** *Assume that we are in the situation of Theorem 8.4. There is a positive integer $c$ such that whenever $k_1$ and $k_2$ are two integers satisfying the congruence*

$$k_1 \equiv k_2 \pmod{(q-1)p^{c[s+1]^{n+r}}},$$

*we have the equality*

$$d_t(k_1) = d_t(k_2)$$

*for all $0 \leq t \leq s$.*

In our current generality, Theorem 8.4 and Corollary 8.10 are probably asymptotically best possible. Of course, they could be improved in various special cases. In the special case that the unit root $\sigma$-module $\psi_{\text{unit}}$ is already overconvergent (for instance, the unit root $\sigma$-module arising from the geometric case of an ordinary family of elliptic curves), both Theorem 8.4 and Corollary 8.10 can be improved. We state this improvement here for comparison purpose. It follows directly from Theorem 8.4 by taking $\psi = \psi_{\text{unit}}$ and thus $r = 1$.

**Corollary 8.11.** *Let $\psi_{\text{unit}}$ be an overconvergent rank one unit root $\sigma$-module. Let $\varphi$ be a fixed overconvergent finite rank $\sigma$-modules. Then, for any integer $k$, we can write*

$$L(\psi_{\text{unit}}^k \otimes \varphi, T) = \frac{f_1(k, T)}{f_2(k, T)},$$

*where $f_1(k, T)$ (resp. $f_2(k, T)$) is a family of uniformly entire functions whose Newton polygons lies above the graph of the function*

$$Q(x) = c_5 x^{1+\frac{1}{n}} - c_6 x$$

*for some positive constants $\{c_5, c_6\}$ independent of $k$. Furthermore, there is a positive integer $c$ such that whenever $k_1$ and $k_2$ are two integers satisfying the congruence*

$$k_1 \equiv k_2 \pmod{(q-1)p^{c[s+1]^{n+1}}},$$

*we have the equality*

$$d_t(k_1) = d_t(k_2)$$

*for all $0 \leq t \leq s$.*

This corollary generalizes the quadratic bound in [15] which corresponds to the special case $n = r = 1$. Corollary 8.11 can be further and greatly improved if we assume the much stronger condition that both $\psi_{\text{unit}}$ and $\varphi$ are overconvergent F-crystals in Berthelot's sense, as the L-function in this case is conjectured to be a rational function. Similarly, it seems plausible to expect that in some cases, Theorem 8.4 and Corollary 8.10 could be improved to the same bounds of Corollary 8.11, if we assume that the ambient $\psi$ and $\varphi$ are overconvergent F-crystals in Berthelot's sense regardless $\psi_{\mathrm{unit}}$ is overconvergent or not as a $\sigma$-module. It would be of great interest to prove such an improvement, which would involve the study of a cohomological formula for the L-function of a nuclear (infinite rank) overconvergent F-crystal.